\documentclass[11pt,letterpaper]{article}
\RequirePackage[backend=biber,bibstyle=authoryear,citestyle=authoryear,
natbib=true,sorting=nyt,sortcites=true,maxnames=10,minnames=10,
hyperref=true,abbreviate=false,date=long,isbn=true,eprint=true,
uniquename=init,giveninits=true]{biblatex}
\usepackage{hyphenat,graphicx}
\usepackage[colorlinks,linktocpage,breaklinks]{hyperref}
\makeatletter
\def\@filecolor{blue}
\def\@linkcolor{blue}
\def\@citecolor{blue}
\def\@urlcolor{blue}
\let\@old@citep\citep
\let\@old@citet\citet
\let\@old@citeauthor\citeauthor
\def\citep{\@old@citep*}
\def\citet{\@old@citet*}
\def\citeauthor{\@old@citeauthor*}
\let\cite\citep
\makeatother
\usepackage[top=1.5in,bottom=1.2in,left=1.15in,right=1.15in,letterpaper]%
  {geometry}
\makeatletter
\def\ps@headings{%
    \let\@mkboth\@gobbletwo
    \def\@oddhead{\hss\scshape\shorttitle\hss\reset@font\rmfamily\thepage}
    \def\@evenhead{\reset@font\rmfamily\thepage\hss\scshape\shortauthors\hss}
    \let\@oddfoot\@empty\let\@evenfoot\@empty}
\@twosidetrue
\makeatother
\usepackage[neveradjust]{paralist}
\makeatletter
\@definecounter{@ADL@savenum}
\newcommand\savenum{\setcounter{@ADL@savenum}%
  {\the\@nameuse{c@\@listctr}}}
\newcommand\resumenum{\setcounter{\@listctr}{\arabic{@ADL@savenum}}}
\makeatother
\usepackage[reqno,centertags]{amsmath}
\usepackage{amssymb,xspace,twoopt,xy,mathrsfs}
\xyoption{arrow}
\xyoption{matrix}
\xyoption{curve}
\usepackage[mathcal]{eucal}
\let\sinfty\infty
\input cyracc.def
\let\subsetneq\subset
\let\subset\subseteq

\newcommand\real{\mathbb{R}}
\newcommand\complex{\mathbb{C}}
\newcommand\field{\mathbb{F}}
\newcommand\realp{\real_{>0}}
\newcommand\realnn{\real_{\ge0}}
\newcommand\integer{\mathbb{Z}}
\newcommand\integerp{\integer_{>0}}
\newcommand\integernn{\integer_{\ge0}}

\newcommand\scirc{\raise1pt\hbox{$\,\scriptstyle\circ\,$}}
\newcommand\sscirc{\hbox{$\,\scriptscriptstyle\circ\,$}}
\newcommand\eqdef{\triangleq}

\newcommand\supscr[2]{#1^{\textup{#2}}}
\newcommand\ol[1]{\overline{#1}}
\newcommand\affdual[1]{\supscr{#1}{$*$,aff}}
\newcommand\disjointunion{\mathop{\overset{\sscirc}{\cup}}}
\newcommand\bigdisjointunion{\mathop{\overset{\sscirc}{\bigcup}}}
\newcommand\what[1]{\widehat{#1}\null}

\newcommand\symmgroup[1]{\mathfrak{S}_{#1}}
\newdir{ >}{{}*!/-3\jot/\dir{>}}
\newcommand\slnorm{\lvert}
\newcommand\srnorm{\rvert}
\newcommand\snorm[1]{\slnorm #1\srnorm}

\newcommand\setdef[2]{\{#1\;|\enspace#2\}}
\newcommand\asetdef[2]{\left\{#1\immediate\vphantom{#2}\;\right|
  \left.\immediate\vphantom{#1}\enspace#2\right\}}
\newcommand\ifam[1]{(#1)}
\newcommand\inprod[2]{\langle#1,#2\rangle}  
\newcommand\natpair[2]{\langle#1;#2\rangle}
\newcommand\ie{i.e.,}
\newcommand\eg{e.g.,}
\newcommand\cf{cf.}

\newcommand\resp{resp.}
\newcommand\mathupper[1]{\textup{#1}}
\renewcommand\d[1]{{\normalfont\textrm{d}}#1}
\newcommand\image{\operatorname{image}}
\newcommand\ev{\operatorname{ev}}
\newcommand\Ev{\operatorname{Ev}}
\newcommand\hcohull{\operatorname{hconv}}
\newcommand\id{\operatorname{id}}

\newcommand\hol{\textup{hol}}

\newcommand\pr{\operatorname{pr}}

\newcommand\rank{\operatorname{rank}}

\newcommand\sM{\mathscr{M}}
\newcommand\sN{\mathscr{N}}
\newcommand\sU{\mathscr{U}}

\let\algdual\stardual 
\let\dual\stardual 
\let\topdual\primedual
\newcommand\tderivatzero[2]{{\textstyle\frac{{\normalfont\mathupper{d}}#1}
  {{\normalfont\mathupper{d}}#2}}\big|_{#2=0}}

\makeatletter
\newcommand\linder{\@ifnextchar[{\@ADL@rlinder}{\@ADL@linder}}
\def\@ADL@rlinder[#1]#2{\boldsymbol{D}^{#1}#2}
\newcommand\@ADL@linder[1]{\boldsymbol{D}#1}
\newcommand\plinder{\@ifnextchar[{\@ADL@rplinder}{\@ADL@plinder}}
\def\@ADL@rplinder[#1]#2#3{\boldsymbol{D}^{#1}_{#2}#3}
\newcommand\@ADL@plinder[2]{\boldsymbol{D}_{#1}#2}
\makeatother
\newcommand\map[3]{#1\colon#2\rightarrow#3}

\newcommand\mapdef[5]{\begin{aligned}
  #1\colon&\begin{aligned}[t]#2\end{aligned}\rightarrow
  \begin{aligned}[t]#3\end{aligned}\\&\begin{aligned}[t]#4\end{aligned}
  \mapsto\begin{aligned}[t]#5\end{aligned}\end{aligned}}
\newcommand\mapschar{\mathupper{C}}
\newcommand\C{\mapschar}

\newcommand\cker[1]{\check{\mathupper{Z}}^{#1}}
\newcommand\cim[1]{\check{\mathupper{B}}^{#1}}
\newcommand\ccohom[1]{\check{\mathupper{H}}^{#1}}
\makeatletter
\renewcommand\hom{\@ifstar{\@ADL@cohom}{\@ADL@hom}}
\newcommand\@ADL@cohom[1]{\mathupper{H}^{#1}}
\newcommand\@ADL@hom[1]{\mathupper{H}_{#1}}
\makeatother
\newcommand\AHom{\mathupper{AHom}}
\newcommand\Hom{\mathupper{Hom}}
\newcommand\Aut{\mathupper{Aut}}
\newcommand\vect[1]{\boldsymbol{#1}}
\newcommand\alg[1]{\mathsf{#1}}
\makeatletter
\newcommand\Symalg{\@ifstar{\@ADL@symalgs}{\@ADL@symalgns}}
\newcommand\@ADL@symalgns{\@ifnextchar[{\@ADL@symalg}{\@ADL@@symalg}}
\def\@ADL@symalg[#1]#2{\@ADL@tsymalgsym^{#1}(#2)}
\newcommand\@ADL@@symalg[1]{\@ADL@tsymalgsym(#1)}
\newcommand\@ADL@symalgs{\@ifnextchar[{\@ADL@ksymalgs}{\@ADL@noksymalgs}}
\def\@ADL@ksymalgs[#1]#2{\@ADL@symalgsym^{#1}(#2)}
\newcommand\@ADL@noksymalgs[1]{\@ADL@symalgsym(#1)}
\newcommand\Symalgsymbol[1]{\def\@ADL@symalgsym{#1}}
\Symalgsymbol{\mathupper{S}}
\makeatother
\newcommand\nbhd[1]{\mathcal{#1}}
\makeatletter
\newcommand\oball{\@ifnextchar[{\@ADL@oballarg}{\@ADL@oballnoarg}}
\def\@ADL@oballarg[#1]#2#3{\mathsf{B}_{#1}(#2,#3)}
\newcommand\@ADL@oballnoarg[2]{\mathsf{B}(#1,#2)}
\newcommand\cball{\@ifnextchar[{\@ADL@cballarg}{\@ADL@cballnoarg}}
\def\@ADL@cballarg[#1]#2#3{\ol{\mathsf{B}}_{#1}(#2,#3)}
\newcommand\@ADL@cballnoarg[2]{\ol{\mathsf{B}}(#1,#2)}
\makeatother
\newcommandtwoopt\sections[3][\sinfty][\null]{\Gamma^{#1}_{#2}(#3)}

\newcommand\man[1]{\mathsf{#1}}
\makeatletter
\newcommand\tb{\@ifnextchar[{\@ADL@tbarg}{\@ADL@tb}}
\def\@ADL@tbarg[#1]#2{\man{T}_{#1}#2}
\newcommand\@ADL@tb[1]{\man{T}#1}
\newcommand\tbproj[1]{\pi_{\tb{#1}}}
\newcommand\ctb{\@ifnextchar[{\@ADL@ctbarg}{\@ADL@ctb}}
\def\@ADL@ctbarg[#1]#2{\man{T}^*_{#1}#2}
\newcommand\@ADL@ctb[1]{\man{T}^*#1}

\newcommand\nb{\@ifnextchar[{\@ADL@nbarg}{\@ADL@nb}}
\def\@ADL@nbarg[#1]#2{\man{N}_{#1}#2}
\newcommand\@ADL@nb[1]{\man{N}#1}
\newcommand\vb{\@ifnextchar[{\@ADL@vbarg}{\@ADL@vb}}
\def\@ADL@vbarg[#1]#2{\man{V}_{#1}#2}
\newcommand\@ADL@vb[1]{\man{V}#1}
\newcommand\cvb{\@ifnextchar[{\@ADL@cvbarg}{\@ADL@cvb}}
\def\@ADL@cvbarg[#1]#2{\man{V}^*_{#1}#2}
\newcommand\@ADL@cvb[1]{\man{V}^*#1}
\newcommand\hb{\@ifnextchar[{\@ADL@hbarg}{\@ADL@hb}}
\def\@ADL@hbarg[#1]#2{\man{H}_{#1}#2}
\newcommand\@ADL@hb[1]{\man{H}#1}
\newcommand\chb{\@ifnextchar[{\@ADL@chbarg}{\@ADL@chb}}
\def\@ADL@chbarg[#1]#2{\man{H}^*_{#1}#2}
\newcommand\@ADL@chb[1]{\man{H}^*#1}
\newcommand\func{\@ifnextchar[{\@ADL@crfuncs}{\@ADL@cinftyfuncs}}
\def\@ADL@crfuncs[#1]#2{\mapschar^{#1}(#2)}
\newcommand\@ADL@cinftyfuncs[1]{\mapschar^\infty(#1)}
\newcommand\mappings{\@ifnextchar[{\@ADL@crmappings}{\@ADL@cinftymappings}}
\def\@ADL@crmappings[#1]#2#3{\mapschar^{#1}(#2;#3)}
\newcommand\@ADL@cinftymappings[2]{\mapschar^\infty(#1;#2)}
\newcommand\Diff{\@ifnextchar[{\@ADL@crdiff}{\@ADL@cinftydiff}}
\def\@ADL@crdiff[#1]#2{\mathupper{Diff}^{#1}(#2)}
\newcommand\@ADL@cinftydiff[2]{\mathupper{Diff}^\sinfty(#1)}
\newcommand\abmappings{\@ifnextchar[{\@ADL@crabmaps}{\@ADL@cinftyabmaps}}
\def\@ADL@crabmaps[#1]#2#3{\mathupper{AB}^{#1}(#2;#3)}
\newcommand\@ADL@cinftyabmaps[2]{\mathupper{AB}^\sinfty(#1;#2)}
\newcommand\vbmappings{\@ifnextchar[{\@ADL@crvbmaps}{\@ADL@cinftyvbmaps}}
\def\@ADL@crvbmaps[#1]#2#3{\mathupper{VB}^{#1}(#2;#3)}
\newcommand\@ADL@cinftyvbmaps[2]{\mathupper{VB}^\sinfty(#1;#2)}
\newcommand\afffunc{\@ifnextchar[{\@ADL@crafffuncs}{\@ADL@cinftyafffuncs}}
\def\@ADL@crafffuncs[#1]#2{\mathupper{Aff}^{#1}(#2)}
\newcommand\@ADL@cinftyafffuncs[1]{\mathupper{Aff}^\sinfty(#1)}
\newcommand\linfunc{\@ifnextchar[{\@ADL@crlinfuncs}{\@ADL@cinftylinfuncs}}
\def\@ADL@crlinfuncs[#1]#2{\mathupper{Lin}^{#1}(#2)}
\newcommand\@ADL@cinftylinfuncs[1]{\mathupper{Lin}^\sinfty(#1)}
\newcommand\gfunc{\@ifnextchar[{\@ADL@crgfuncs}{\@ADL@cinftygfuncs}}
\def\@ADL@crgfuncs[#1]#2#3{\mathscr{C}^{#1}_{#2,#3}}
\newcommand\@ADL@cinftygfuncs[2]{\mathscr{C}^\infty_{#1,#2}}
\newcommand\sfunc{\@ifnextchar[{\@ADL@crsfuncs}{\@ADL@cinftysfuncs}}
\def\@ADL@crsfuncs[#1]#2{\mathscr{C}^{#1}_{#2}}
\newcommand\@ADL@cinftysfuncs[1]{\mathscr{C}^\infty_{#1}}
\newcommand\gsections{\@ifnextchar[{\@ADL@rgsections}{\@ADL@gsections}}
\def\@ADL@rgsections[#1]#2#3{\mathscr{G}^{#1}_{#2,#3}}
\newcommand\@ADL@gsections[2]{\mathscr{G}^\infty_{#1,#2}}
\newcommand\ssections{\@ifnextchar[{\@ADL@rssections}{\@ADL@ssections}}
\def\@ADL@rssections[#1]#2{\mathscr{G}^{#1}_{#2}}
\newcommand\@ADL@ssections[1]{\mathscr{G}^\infty_{#1}}
\newcommand\tf{\@ifnextchar[{\@ADL@tfarg}{\@ADL@tf}}
\def\@ADL@tfarg[#1]#2{\@ifnextchar[{\@ADL@@tfargr{#1}{#2}}
  {\@ADL@@tfarg{#1}{#2}}}
\def\@ADL@@tfargr#1#2[#3]{T^{#3}_{#1}#2}
\newcommand\@ADL@@tfarg[2]{T_{#1}#2}
\newcommand\@ADL@tf[1]{\@ifnextchar[{\@ADL@@tfr{#1}}{\@ADL@@tf{#1}}}
\def\@ADL@@tfr#1[#2]{T^{#2}#1}
\newcommand\@ADL@@tf[1]{T#1}
\newcommand\ctf{\@ifnextchar[{\@ADL@ctfarg}{\@ADL@ctf}}
\def\@ADL@ctfarg[#1]#2{T^*_{#1}#2}
\newcommand\@ADL@ctf[1]{T^*#1}
\newcommand\lieder[2]{\def\@tempa{#2}\ifx\@tempa\@empty%
  \boldsymbol{\mathscr{L}}_{#1}\else\boldsymbol{\mathscr{L}}_{#1}#2\fi}
\newcommand\hl{\@ifstar{\@hl@star}{\@hl@nostar}}
\newcommand\@hl@star[1]{\supscr{#1}{h,$*$}}
\newcommand\@hl@nostar[1]{\supscr{#1}{h}}
\makeatother
\newcommand\ve[1]{\supscr{#1}{e}}

\makeatletter
\def\flow{\@ifstar{\@flow@star}{\@flow@nostar}}
\def\@flow@star#1{\Phi^{#1}}
\def\@flow@nostar#1#2{\@ifnextchar[{\@tflow{#1}{#2}}{\@flow{#1}{#2}}}
\def\@tflow#1#2[#3]{\Phi^{#1}_{#2,#3}}
\def\@flow#1#2{\Phi^{#1}_{#2}}
\def\jet{\@ifnextchar[{\@ADL@jetbase}{\@ADL@jetnobase}}
\def\@ADL@jetbase[#1]#2#3{\man{J}_{#1}^{#2}#3}
\def\@ADL@jetnobase#1#2{\man{J}^{#1}#2}
\def\jetalg{\@ifnextchar[{\@ADL@jetalgbase}{\@ADL@jetalgnobase}}
\def\@ADL@jetalgbase[#1]#2#3{\man{T}^{*#2}_{#1}#3}
\def\@ADL@jetalgnobase#1#2{\man{T}^{*#1}#2}
\makeatother
\newcommandtwoopt\DO[3][\sinfty][\null]{\mathupper{DO}^{#1}_{#2}(#3)}
\newcommandtwoopt\FADO[3][\sinfty][\null]{\mathupper{FADO}^{#1}_{#2}(#3)}
\newcommandtwoopt\FLDO[3][\sinfty][\null]{\mathupper{FLDO}^{#1}_{#2}(#3)}
\newcommandtwoopt\ADO[3][\sinfty][\null]{\mathupper{ADO}^{#1}_{#2}(#3)}
\newcommandtwoopt\LDO[3][\sinfty][\null]{\mathupper{LDO}^{#1}_{#2}(#3)}
\makeatletter
\def\interval{\@ifnextchar({\@ADL@openleftint}{\@ADL@closedleftint}}
\def\@ADL@openleftint(#1,#2{(#1,#2%
  \@ifnextchar){\@ADL@openrightint}{\@ADL@closedrightint}}
\def\@ADL@closedleftint[#1,#2{[#1,#2%
  \@ifnextchar){\@ADL@openrightint}{\@ADL@closedrightint}}
\def\@ADL@openrightint){)}
\def\@ADL@closedrightint]{]}
\makeatother
\usepackage{theorem}
\newtheorem{theorem}{Theorem}[section]
\newtheorem{proposition}[theorem]{Proposition}
\newtheorem{lemma}[theorem]{Lemma}
\newtheorem{corollary}[theorem]{Corollary}
\newtheorem{prooflemma}{Lemma}[theorem]

{\theorembodyfont{\normalfont\rmfamily}
\newtheorem{definition}[theorem]{Definition}

\newtheorem{examples}[theorem]{Examples}}
\makeatletter
\newcommand\pushright{\protect\@ADL@pushright}
\newcommand\@ADL@pushright[1]{{\ifvmode\null\hfill{#1}\par\else\ifmmode%
  \@ADLmaths@pushright{\hbox{#1}}\else\ifinner\@ADLhbox@pushright{#1}%
  \else\@ADLparag@pushright{#1}\fi\fi\fi}}
\newcommand\@ADLmaths@pushright[1]{{\ifinner\@ADLhbox@pushright{#1}\else%
  \tag*{$#1$}\fi}}
\newcommand\@ADLparag@pushright[1]{{\parfillskip=0pt\widowpenalty=10000%
  \displaywidowpenalty=10000\finalhyphendemerits=0\@ADLhbox@pushright#1\par}}
\newcommand\@ADLhbox@pushright{\unskip\nobreak\hfil\penalty50\hskip.2em%
  \null\hfill\hfill}
\newenvironment{proof}{\trivlist\item[\hskip\labelsep\textit{Proof:}\/]%
  \@ADLsave@set@qed\xspace\normalfont\rmfamily}
  {\qed\@ADLrestore@qed\endtrivlist}
\newif\if@ADL@qed\@ADL@qedfalse
\newcommand\qed{\protect\@ADL@qed{$\blacksquare$}}
\newcommand\@ADL@qed[1]{\if@ADL@qed\global\@ADL@qedfalse%
  \pushright{#1}\else\ifhmode\ifinner\else\par\fi\fi\fi}
\newcommand\@ADLrestore@qed{\global\let\if@ADL@qed\@ADLsaved@ifqed}
\newcommand\@ADLsave@set@qed{\let\@ADLsaved@ifqed
  \if@ADL@qed\global\@ADL@qedtrue}

\newenvironment{subproof}{\trivlist\item[\hskip\labelsep\textit{Proof:}\/]%
  \@ADLsave@set@subqed\normalfont\rmfamily}
  {\subqed\@ADLrestore@subqed\endtrivlist}
\newif\if@ADL@subqed\@ADL@subqedfalse
\newcommand\subqed{\protect\@ADL@subqed{$\blacktriangledown$}}
\newcommand\@ADL@subqed[1]{\if@ADL@subqed\global\@ADL@subqedfalse%
  \pushright{#1}\else\ifhmode\ifinner\else\par\fi\fi\fi}
\newcommand\@ADLrestore@subqed{\global\let\if@ADL@subqed\@ADLsaved@ifsubqed}
\newcommand\@ADLsave@set@subqed{\let\@ADLsaved@ifsubqed
  \if@ADL@subqed\global\@ADL@subqedtrue}

\newif\if@ADL@oprocend\@ADL@oprocendfalse
\newcommand\oprocend{\@ADLsave@set@oprocend
  \protect\@ADL@oprocend{$\bullet$}\@ADLrestore@oprocend}
\newcommand\@ADL@oprocend[1]{\if@ADL@oprocend\global\@ADL@oprocendfalse%
  \pushright{#1}\else\ifhmode\ifinner\else\par\fi\fi\fi}
\newcommand\@ADLrestore@oprocend{\global
  \let\if@ADL@oprocend\@ADLsaved@ifoprocend}
\newcommand\@ADLsave@set@oprocend{\let\@ADLsaved@ifoprocend\if@ADL@oprocend%
  \global\@ADL@oprocendtrue}

\makeatother
\newenvironment{keywords}{\quote\small\textbf{Keywords.}}{\endquote}
\newenvironment{AMS}{\quote\small\textbf{AMS Subject Classifications (2010).}}
   {\endquote}
\newcommand\defn[1]{{\normalfont\bfseries\emph{\mathversion{bold}#1}}}
\makeatletter
\@namedef{procnonum*}{\trivlist \@topsep \theorempreskipamount
  \@topsepadd \theorempostskipamount
  \@ifnextchar[{\@ADL@xprocnonumstar}{\@ADL@yprocnonumstar}}
\@namedef{endprocnonum*}{\endtrivlist}
\def\@ADL@xprocnonumstar[#1]{\item[\hskip \labelsep{\theorem@headerfont #1}]
  \normalfont\rmfamily}
\def\@ADL@yprocnonumstar{\item[] \normalfont\rmfamily}
\makeatother
\numberwithin{equation}{section}

\newcommand\enumdblref[2]{\mbox{\ref{#1}--\ref{#2}}}
\title{Gelfand duality for manifolds, and vector and other bundles\thanks{Research supported in part by a grant from the Natural Sciences and Engineering Research Council of Canada}}
\author{Andrew D.\ Lewis\thanks{Professor, Department of Mathematics and
Statistics, Queen's University, Kingston, ON K7L 3N6, Canada,
email:~\texttt{andrew.lewis@queensu.ca}}}
\date{2020/09/20}
\newcommand\shorttitle{Gelfand duality for manifolds, and vector and other bundles}
\newcommand\shortauthors{A.\ D.\ Lewis}
\begin{document}
\maketitle

\begin{abstract}
In general terms, Gelfand duality refers to a correspondence between a
geometric, topological, or analytical category, and an algebraic category.
For example, in smooth differential geometry, Gelfand duality refers to the
topological embedding of a smooth manifold in the topological dual of its
algebra of smooth functions.  This is generalised here in two directions.
First, the topological embeddings for manifolds are generalised to the cases
of real analytic and Stein manifolds, using a unified cohomological argument.
Second, this type of duality is extended to vector bundles, affine bundles,
and jet bundles by using suitable classes of functions, the topological duals
in which the embeddings take their values.
\end{abstract}

\section{Introduction}\label{sec:intro}

In smooth differential geometry, the embedding of a smooth manifold $\man{M}$
in the topological dual $\topdual{\func[\infty]{\man{M}}}$ of its algebra of
smooth functions is well known~\cite[\eg][Theorem~7.2]{JN:03}\@.  That the
result is true appears to originate as~\cite[Problem~1-C]{JWM/JDS:74}\@.
This is also known in the case of Stein manifolds, and seems due
to~\citet[Theorem~2.6]{HR:63}\@.  For real analytic manifolds, the first
proof of which we are aware is given in the PhD thesis of
\citet[Theorem~3.4.4]{SJ:16}\@.  In all cases, the idea is a refrain from
algebraic geometry: one wishes to understand correspondences between a space
and the space of natural functions defined on the space.  This general idea
is known as Gelfand duality, and can be regarded as providing a full and
faithful functor between the geometric category (say, smooth manifolds and
mappings) and the algebraic category (say, the opposite category\footnote{It
has to be the opposite category because of the fact that the pull-back of a
composition is the reversed composition of the pull-backs.} of the category
of $\real$-algebras).

In this paper we are concerned with the geometry/algebra correspondence for
bundles with algebraic structure, specifically vector, affine, and jet
bundles.  For vector bundles, certain algebraic correspondences are standard.
\begin{compactenum}
\item \emph{The ``locally free, locally finitely generated''
correspondence}\@: A commonly called upon correspondence between a geometric
object and a (sort of) algebraic object is the categorical equivalence
between vector bundles and locally free, locally finitely generated sheaves
of modules over the sheaf of rings of functions.  This is explained in the
context of smooth geometry by \citet[\S2.2]{SR:05} and quickly by
\citet[\S1.4.2]{HG/RR:84} in the holomorphic case.
\item \emph{The ``projective module'' correspondence}\@: A related
correspondence that is brought up in the same vein is the correspondence
between the modules of sections of vector bundles and finitely generated
projective modules over the ring of functions.  This is known as the
``Serre\textendash{}Swan Theorem'' as it is proved for algebraic vector
bundles over affine varieties by \citet{JPS:55} and for vector bundles over
compact Hausdorff topological spaces by \citet{RGS:62}\@.  The version for
smooth vector bundles came into being at some point, and is given by
\cite[Theorem~11.32]{JN:03}\@.
\end{compactenum}

These well-known geometric/algebraic correspondences for vector bundles are
not without their limitations.  The correspondence with locally free, locally
finitely generated sheaves is quite perfect; indeed, it is rather close to a
tautology once one understands the words involved.  On the other hand, the
correspondence with finitely generated projective modules is not quite
tautological.  However, it is not uniquely defined, in the sense that the
module of sections can be a summand of a module in many different ways.
Also, this projective module characterisation is of a different character
than the standard Gelfand correspondence for manifolds.  Indeed, the two
correspondences seem a bit orthogonal.  Moreover, both of these
correspondences become complicated when the one talks about vector bundles
over different base spaces, as the base ring changes for the modules under
consideration.

In the paper we address three questions that arise from the preceding
discussion.  
\begin{compactenum}
\item Can the Gelfand duality for manifolds be unified across regularity
classes?
\item Does the Serre\textendash{}Swan Theorem hold for vector bundles with
regularity other than what is mentioned above?
\item Is there a full and faithful functor for vector (or other) bundles that
more closely resembles the Gelfand duality for manifolds?
\end{compactenum}

As to the first question, we give a proof of Gelfand duality for smooth, real
analytic, and Stein manifolds that is ``the same'' for all cases.  It relies
on reducing a crucial part of the proof to an argument using the vanishing of
sheaf cohomology in the three cases.  The results are presented in
Section~\ref{sec:Membedding}\@, with the main result being Theorem~\ref{the:Membedding}\@.

Concerning the second question, we prove as Theorem~\ref{the:swan} the
Serre\textendash{}Swan Theorem for smooth, real analytic, and Stein base
spaces, again using a unified argument.

The answering of these first two questions can be seen as wrapping up some
loose ends, and tightening up the presentation of existing results.  However,
the line of the third question seems unaddressed in the existing literature.
In Sections~\ref{sec:Aembedding} and~\ref{sec:Jembedding} we answer the
question by providing, for a few classes of bundles\textemdash{}vector
bundles, affine bundles, and jet bundles\textemdash{}a version of Gelfand
duality for these spaces.  The geometry/algebra correspondence we give is
closely integrated with the standard Gelfand duality for manifolds, which
distinguishes it from the existing correspondences for vector bundles.
Another feature of our approach to the geometry/algebra correspondence for
bundles is that it utilises topological properties of the function spaces
involved to make the correspondences homeomorphisms. The key idea is the
determination of (1)~an appropriate space of functions to play the r\^ole of
the algebra of all functions with desired regularity and (2)~the appropriate
set of morphisms that makes the correspondence functorial.

\subsection*{Notation}

When $A$ is a subset of a set $X$\@, we write $A\subset X$\@.  If we wish to
exclude the possibility that $A=X$\@, we write $A\subsetneq X$\@.  The
identity map on a set $X$ is denoted by $\id_X$\@.

By $\integer$ we denote the set of integers.  We use the notation $\integerp$
and $\integernn$ to denote the subsets of positive and nonnegative integers.
By $\real$ we denote the sets of real numbers.  By $\realp$ we denote the
subset of positive real numbers.  By $\complex$ we denote the set of complex
numbers.  We shall work simultaneously with real and complex numbers, and so
denote $\field\in\{\real,\complex\}$ in these cases.  We denote by $\field^n$
the $n$-fold Cartesian product of $\field$\@.

If $\alg{R}$ is a ring (a commutative ring with unit) and if $\alg{U}$ and
$\alg{V}$ and $\alg{R}$-modules, we denote by
$\Hom_{\alg{R}}(\alg{U};\alg{V})$ the set of module homomorphisms from
$\alg{U}$ to $\alg{V}$\@.  We denote by
$\dual{\alg{V}}=\Hom_{\alg{R}}(\alg{V};\alg{R})$ the algebraic dual.  If
$v\in\alg{V}$ and $\alpha\in\dual{\alg{V}}$\@, we will denote the evaluation
of $\alpha$ on $v$ at various points by $\alpha(v)$\@, $\alpha\cdot v$\@, or
$\natpair{\alpha}{v}$\@, whichever seems most pleasing to us at the moment.

By $\symmgroup{k}$ we denote the permutation group of $\{1,\dots,k\}$\@.  For
$k,l\in\integernn$\@, we denote by $\symmgroup{k,l}$ the subset of
$\symmgroup{k+l}$ consisting of permutations $\sigma$ satisfying
\begin{equation*}
\sigma(1)<\dots<\sigma(k),\quad\sigma(k+1)<\dots<\sigma(k+l).
\end{equation*}

By $\Symalg*[k]{\alg{V}}$ we denote the $k$-fold symmetric tensor product of
$\alg{V}$ with itself, and we think of this as a subset of the $k$-fold
tensor product.  For $A\in\Symalg*[k]{\alg{V}}$ and
$B\in\Symalg*[l]{\alg{V}}$\@, we define the symmetric tensor product of $A$
and $B$ to be
\begin{equation*}
A\odot B=\sum_{\sigma\in\symmgroup{k,l}}\sigma(A\otimes B).
\end{equation*}

We shall adopt the notation and conventions of smooth differential geometry
of~\cite{RA/JEM/TSR:88}\@.  We shall also make use of real analytic
differential geometry.  There are no useful textbook references dedicated to
real analytic differential geometry, but the book of \citet{KC/YE:12}
contains much of what we shall need.  For complex geometry, we refer
to~\cite{ROW:08}\@. Throughout the paper, manifolds are connected, second
countable, Hausdorff manifolds.  The assumption of connectedness can be
dispensed with but is convenient as it allows one to not have to worry about
manifolds with components of different dimensions and vector bundles with
fibres of different dimensions.

We shall work with regularity classes $r\in\{\infty,\omega,\hol\}$\@,
``$\infty$'' meaning smooth, ``$\omega$'' meaning real analytic, and
``$\hol$'' meaning holomorphic.  We shall use $\field=\real$ when working in
the smooth and real analytic settings, and use $\field=\complex$ when working
in the holomorphic setting.  In the holomorphic case, we work with Stein
manifolds and with vector bundles over Stein manifolds.  When we are being
careful, as in stating theorems, we will be sure to state this clearly.
However, in discussions, we will sometimes make statements about holomorphic
geometry that are only true for Stein manifolds or for vector bundles over
Stein manifolds, while not explicitly mentioning that the Stein assumption is
being made.  For this reason, it is probably best to always assume the Stein
assumption is being made in the background when it is not being made in the
foreground.

We denote by $\mappings[r]{\man{M}}{\man{N}}$ the set of mappings from a
manifold $\man{M}$ to a manifold $\man{N}$ of class $\C^r$\@.  By
$\Diff[r]{\man{M}}$ we denote the set of $\C^r$-diffeomorphisms of
$\man{M}$\@.  When $\man{N}=\field$\@, we denote by
$\func[r]{\man{M}}=\mappings[r]{\man{M}}{\field}$ the set of scalar-valued
functions of class $\C^r$\@.  We denote by $\mathsf{1}_{\man{M}}$ the
constant function with value $1$ on a manifold $\man{M}$\@.  By $\d{f}$ we
denote the differential of $f$\@.

By $\map{\tbproj{\man{M}}}{\tb{\man{M}}}{\man{M}}$ we denote the tangent
bundle of $\man{M}$ (the holomorphic tangent bundle in the holomorphic
case).  If $\Phi\in\mappings[r]{\man{M}}{\man{N}}$\@, we denote by
$\map{\tf{\Phi}}{\tb{\man{M}}}{\tb{\man{N}}}$ the derivative of $\Phi$\@.  By
$\tf[x]{\Phi}$ we denote the restriction of $\tf{\Phi}$ to
$\tb[x]{\man{M}}$\@.

Let $\map{\pi}{\man{E}}{\man{M}}$ be a vector bundle of class $\C^r$\@.  We
shall sometimes denote the fibre over $x\in\man{M}$ by $\man{E}_x$\@, noting
that this has the structure of an $\field$-vector space.  By
$\field_{\man{M}}=\man{M}\times\field$\@, we denote the trivial line bundle.
If $A\subset\man{M}$\@, we denote by $\man{E}|A=\pi^{-1}(A)$\@.  If
$\man{S}\subset\man{M}$ is a submanifold, then $\man{E}|\man{S}$ is a vector
bundle over $\man{S}$\@.  By $\sections[r]{\man{E}}$ we denote the set of
sections of $\man{E}$ of class $\C^r$\@.  If $\map{\pi}{\man{E}}{\man{M}}$
and $\map{\theta}{\man{F}}{\man{N}}$ are $\C^r$-vector bundles, a
$\C^r$-vector bundle mapping from $\man{E}$ to $\man{F}$ is a pair
$(\Phi,\Phi_0)$ of $\C^r$-mappings making the diagram
\begin{equation*}
\xymatrix{{\man{E}}\ar[r]^{\Phi}\ar[d]_{\pi}&{\man{F}}\ar[d]^{\theta}\\
{\man{M}}\ar[r]_{\Phi_0}&{\man{N}}}
\end{equation*}
commute, and such that
$\Phi_x\in\Hom_{\field}(\man{E}_x;\man{F}_{\Phi_0(x)})$\@, where
$\Phi_x=\Phi|\man{E}_x$\@.  We denote by $\vbmappings[r]{\man{E}}{\man{F}}$
the set of vector bundle mappings from $\man{E}$ to $\man{F}$\@.

By $\sfunc[r]{\man{M}}$ we denote the sheaf of $\C^r$-function on $\man{M}$
and by $\ssections[r]{\man{E}}$ we denote the sheaf of sections of
$\man{E}$\@, thought of as an $\sfunc[r]{\man{M}}$-module.

We shall often make use of the fact that, for the manifolds we consider,
there are always globally defined coordinate functions.
\begin{lemma}\label{lem:global-coordinates}
Let\/ $r\in\{\infty,\omega,\hol\}$ and let\/
$\field\in\{\real,\complex\}$\@, as appropriate.  If\/ $\man{M}$ is a\/
$\C^r$-manifold, Stein when\/ $r=\hol$\@, then, for any\/ $x\in\man{M}$\@,
there exists a chart\/ $(\nbhd{U},\phi)$ for\/ $\man{M}$ whose coordinate
functions\/ $\chi^1,\dots,\chi^n$ are restrictions to\/ $\nbhd{U}$ of
globally defined functions of class\/ $\C^r$\@.
\begin{proof}
The hypotheses ensure that there is a proper $\C^r$-embedding
$\map{\iota_{\man{M}}}{\man{M}}{\field^N}$ for some suitable
$N\in\integerp$\@; this is a result of \citet{HW:36} in the smooth case,
\citet{HG:58} in the real analytic case, and \citet{RR:54} in the holomorphic
case.  Define $\chi^1,\dots,\chi^N\in\func[r]{\man{M}}$ by
\begin{equation*}
\iota_{\man{M}}(x)=(\chi^1(x),\dots,\chi^N(x)),\qquad x\in\man{M}.
\end{equation*}
Now, for $x\in\man{M}$\@, $\tf[x]{\iota_{\man{M}}}$ is injective, and so
there exists $j_1,\dots,j_n\in\{1,\dots,N\}$ such that
$\ifam{\d{\chi^{j_1}}(x),\dots,\d{\chi^{j_n}}(x)}$ is a basis for
$\ctb[x]{\man{M}}$\@.  In some neighbourhood $\nbhd{U}$ of $x$\@, we will
have linear independence of
$\ifam{\d{\chi^{j_1}}(x),\dots,\d{\chi^{j_n}}(x)}$\@, and so
$\chi^{j_1},\dots,\chi^{j_n}$ are coordinate functions on $\nbhd{U}$\@.
\end{proof}
\end{lemma}

Finally, we mention that for the topological assertions in our main results,
we make use of topologies on spaces of $\C^r$-vector bundles.  We do not go
into detail about what these topologies are, as a detailed understanding of
this is not material to the important points we are making in this paper.
However, there may be some applications of our results here that will benefit
from a detailed understanding of the topologies involved.  We refer
to~\cite{ADL:20a} for details in the real analytic case.  Also in that work
some words are said about how to simplify the proofs in the real analytic
case to the smooth case.  The holomorphic case is even easier, as the
topology in this case is the topology of uniform convergence on compact
sets.  Alternatively, one can use the smooth topology on holomorphic mappings
since holomorphic mappings form a closed subset of smooth mappings~\cite[Theorem~II.8.2]{AK/PWM:97}\@.

\section{Functions on vector, affine, and jet bundles}

While vector bundles are most commonly encountered in differential geometry,
for what we do in this paper it is most natural to work with affine bundles.
The reason for this is that it is affine functions (not linear functions)
that we will use to characterise Gelfand duality for vector bundles.
Therefore, we prefer to work with affine bundles, where affine functions are
most naturally defined.  Also, the affine structure of jet bundles makes it
possible to extend our notions of Gelfand duality from affine bundles to jet bundles.

\subsection{Affine bundles}

We assume the reader is familiar with the notion of an affine space modelled
on a vector space~\cite[Chapter~2]{MB:87}\@.  If $\alg{A}$ is an affine space
modelled on a vector space $\alg{V}$\@, then $\affdual{\alg{A}}$ denotes the
\defn{affine dual} of $\alg{A}$\@, by which we mean the set of affine maps
from $\alg{A}$ to the field over which $\alg{V}$ is defined.  We shall adopt
the notation for linear duals, and write $\natpair{\lambda}{a}=\lambda(a)$
for the evaluation of $\lambda\in\affdual{\alg{A}}$ on $a\in\alg{A}$\@.

An affine bundle is the extension of the idea of an affine space to
differential geometry, in the same way as a vector bundle is an extension to
differential geometry of a vector space.
\begin{definition}
Let $r\in\{\infty,\omega,\textup{hol}\}$ and let $\field=\real$ if
$r\in\{\infty,\omega\}$ and let $\field=\complex$ if $r=\textup{hol}$\@.  Let
$\map{\pi}{\man{E}}{\man{M}}$ be a $\C^r$-vector bundle.  A
\defn{$\C^r$-affine bundle} over $\man{M}$ modelled on $\man{E}$ is a
$\C^r$-fibre bundle $\map{\beta}{\man{B}}{\man{M}}$ with the following
structure:
\begin{compactenum}[(i)]
\item there exists a $\C^r$-fibre bundle mapping
$\map{\alpha}{\man{E}\times_{\man{M}}\man{B}}{\man{B}}$ such that
$a+e\eqdef\alpha(e,a)$ makes $\man{B}_x$ into an affine space modelled on
$\man{E}_x$ for each $x\in\man{M}$\@;
\item for each\/ $x\in\man{M}$\@, there exists a $\C^r$-local trivialisation
$\map{\tau}{\beta^{-1}(\nbhd{U})}{\nbhd{U}\times\field^k}$ for which
$\map{\pr_2\scirc(\tau|\man{B}_y)}{\man{B}_y}{\field^k}$ is an isomorphism of
affine spaces for each $y\in\nbhd{U}$\@.\oprocend
\end{compactenum}
\end{definition}

Let us flesh out the meaning of the second condition concerning local
trivialisations.  Suppose that $(\nbhd{U},\phi)$ is a vector bundle chart for
the model vector bundle $\map{\pi}{\man{E}}{\man{M}}$ associated with the
affine bundle $\map{\beta}{\man{B}}{\man{M}}$\@.  Then, possibly after
shrinking $\nbhd{U}$\@, there is a local trivialisation
$\map{\tau}{\beta^{-1}(\nbhd{U}_0)}{\nbhd{U}_0\times\field^k}$
satisfying the second condition in the definition.  We also have a local
trivialisation
$\map{\lambda}{\pi^{-1}(\nbhd{U}_0)}{\nbhd{U}_0\times\field^k}$
given by the vector bundle chart that is a vector bundle mapping.  This means
that the representation in these local trivialisations of the mapping
$\alpha$ providing the affine structure is
\begin{equation*}
(\vect{x},(\vect{e},\vect{a}))\mapsto(\vect{x},\vect{a}+\vect{e}).
\end{equation*}
Thus, locally, the affine bundle looks like the product of an open set with
the affine space~$\field^k$\@.

\begin{definition}
Let $r\in\{\infty,\omega,\textup{hol}\}$ and let $\field=\real$ if
$r\in\{\infty,\omega\}$ and let $\field=\complex$ if $r=\textup{hol}$\@.  If
$\map{\beta_1}{\man{B}_1}{\man{M}_1}$ and
$\map{\beta_2}{\man{B}_2}{\man{M}_2}$ are $\C^r$-affine bundles, then a
\defn{$\C^r$-affine bundle map} between these affine bundles is a $\C^r$-map
$\map{\Phi}{\man{B}_1}{\man{B}_2}$ for which there exists a $\C^r$-map
$\map{\Phi_0}{\man{M}_1}{\man{M}_2}$ such that the diagram
\begin{equation*}
\xymatrix{{\man{B}_1}\ar[r]^{\Phi}\ar[d]_{\beta_1}&{\man{B}_2}\ar[d]^{\beta_2}\\
{\man{M}_1}\ar[r]_{\Phi_0}&{\man{M}_2}}
\end{equation*}
commutes and with the property that
$\map{\Phi|\man{B}_{1,x}}{\man{B}_{1,x}}{\man{B}_{2,\Phi_0(x)}}$ is an affine
map.  If $\Phi$ is a $\C^r$-diffeomorphism we say it is an \defn{affine
bundle isomorphism}\@.

We denote by $\abmappings[r]{\man{B}_1}{\man{B}_2}$ the set of $\C^r$-affine
bundle mappings from $\man{B}_1$ to $\man{B}_2$\@.\oprocend
\end{definition}

We let $\sections[r]{\man{B}}$ denote the set of sections of an affine bundle
$\man{B}$\@.  Being fibre bundles, affine bundles are entitled to the
possession of local sections.  Also, just being fibre bundles, they are not
\emph{a priori} entitled the possession of global sections.  However, one
feels that they are really a lot like vector bundles, and so should possess
as many global sections as their model vector bundles.  Indeed, if an affine
bundle $\map{\beta}{\man{B}}{\man{M}}$ possesses \emph{one} section
$\sigma$\@, then $\sigma+\xi$ is also a section for any section $\xi$ of the
model vector bundle $\map{\pi}{\man{E}}{\man{M}}$\@.  Thus the question of
the character of the set of sections of an affine bundle really boils down to
the existence on one section.
\begin{proposition}\label{prop:affbunsec}
Let\/ $r\in\{\infty,\omega,\textup{hol}\}$ and let\/ $\field=\real$ if\/
$r\in\{\infty,\omega\}$ and let\/ $\field=\complex$ if\/ $r=\textup{hol}$\@.
Let\/ $\map{\beta}{\man{B}}{\man{M}}$ be a\/ $\C^r$-affine bundle modelled on
the\/ $\C^r$-vector bundle\/ $\map{\pi}{\man{E}}{\man{M}}$ and suppose that\/
$\man{M}$ is Stein if\/ $r=\hol$\@.  Then\/
$\sections[r]{\man{B}}\not=\emptyset$\@.
\begin{proof}
To begin, we claim that there is an open cover $\ifam{\nbhd{U}_a}_{a\in A}$
such that the sheaves $\ssections[r]{\man{E}}|(\cap_{a\in F}\nbhd{U}_a)$\@,
$F\subset A$ finite, are acyclic.  First of all, in the smooth and real
analytic cases, this is true for \emph{every} open cover.  In the smooth
case, this follows from the vanishing of the sheaf cohomology for sheaves of
modules over $\sfunc[\infty]{\man{M}}$ (\cite[Proposition~3.11]{ROW:08}\@,
along with \cite[Examples~3.4(d,e)]{ROW:08}
and~\cite[Proposition~3.5]{ROW:08}).  In the real analytic case, we note that
$\ssections[\omega]{\man{E}}$ is coherent by the Oka Coherence Theorem (see
\cite[Theorem~2.5.2]{HG/RR:84} in the holomorphic case; the same proof works
in the real analytic case).  Thus the assertion follows from Cartan's
Theorem~B~\cite[Proposition~6]{HC:57}\@.  In the holomorphic case, we can
assume that the open sets are Stein (\eg~we can take the open sets
$\nbhd{U}_a$\@, $a\in A$\@, to be preimages under an holomorphic chart of a
polydisk in $\complex^n$).  Then, since finite intersections of Stein open
sets are Stein (by \cite[Proposition~I.6.20(c)]{JPD:12}), it follows from
Cartan's Theorem~B that $\ssections[r]{\man{E}}|(\cap_{a\in F}\nbhd{U}_a)$ is
acyclic for every finite $F\subset A$\@.

About any $x\in\man{M}$ there is a neighbourhood $\nbhd{U}$ so that
$\man{E}|\nbhd{U}\simeq\nbhd{U}\times\field^k$\@,
$\man{B}|\nbhd{U}\simeq\nbhd{U}\times\field^k$\@, and the affine structure on
the fibres is the standard one.  Therefore, there are local sections
$\sigma_1,\dots,\sigma_{k+1}\in\sections[r][\nbhd{U}]{\man{B}}$ such that any
local section of $\man{B}|\nbhd{U}$ is an affine combination of these,~\ie
\begin{equation*}
\sigma\in\sections[r][\nbhd{U}]{\man{B}}\enspace\implies\enspace
\sigma(x)=\sum_{a=1}^{k+1}f^a(x)\sigma_a(x),\quad
f^1,\dots,f^{k+1}\in\func[r]{\nbhd{U}},\ \sum_{a=1}^{k+1}f^a(x)=1,\
x\in\nbhd{U}.
\end{equation*}
Therefore, there exists an open cover $\sU=\ifam{\nbhd{U}_a}_{a\in A}$ for
$\man{M}$ such that, for each $a\in A$\@, we have local generators
$\ifam{\sigma_{ai}}_{i\in I}$ for the affine bundle $\man{B}$ (as above) on
$\nbhd{U}_a$\@.  In the holomorphic case, we assume that the open sets
$\nbhd{U}_a$ are Stein, as can be done without loss of generality.  The index
set $I$ can be taken to be the same for all open sets by our assumption that
$\man{M}$ is connected.  For $a\in A$\@, fix $i_0\in I$ and denote
$\sigma_{a0}=\sigma_{ai_0}$\@.  For $x\in\nbhd{U}_a$\@, we have
$\man{B}_x=\sigma_{a0}(x)+\man{E}_x$\@.  If
$\nbhd{U}_a\cap\nbhd{U}_b\not=\emptyset$\@, then we have
$\sigma_{a0}(x)-\sigma_{b0}(x)\in\man{E}_x$ for
$x\in\nbhd{U}_a\cap\nbhd{U}_b$\@.  Said otherwise,
\begin{equation*}
\sigma_{a0}|\nbhd{U}_a\cap\nbhd{U}_b-\sigma_{b0}|\nbhd{U}_a\cap\nbhd{U}_b\in
\ssections[r]{\man{E}}(\nbhd{U}_a\cap\nbhd{U}_b).
\end{equation*}
Denote $\xi_{ab}\in\ssections[r]{\man{E}}(\nbhd{U}_a\cap\nbhd{U}_b)$ by
\begin{equation*}
\xi_{ab}=\sigma_{a0}|\nbhd{U}_a\cap\nbhd{U}_b-
\sigma_{b0}|\nbhd{U}_a\cap\nbhd{U}_b,
\end{equation*}
and note that
\begin{equation*}
\xi_{ac}|\nbhd{U}_a\cap\nbhd{U}_b\cap\nbhd{U}_c=
\xi_{ab}|\nbhd{U}_a\cap\nbhd{U}_b\cap\nbhd{U}_c+
\xi_{bc}|\nbhd{U}_a\cap\nbhd{U}_b\cap\nbhd{U}_c.
\end{equation*}
Thus we have a \v{C}ech $1$-cocycle
$\ifam{\xi_{ab}}_{a,b\in A}\in\cker{1}(\sU;\ssections[r]{\man{E}})$\@.  As we
have seen in the first paragraph of the proof,
$\hom*{1}(\man{M};\ssections[r]{\man{E}})=0$\@.  By Leray's
Theorem~\cite[Theorem~5.3]{SR:05}\@, the \v{C}ech cohomology
$\ccohom{1}(\sU;\ssections[r]{\man{E}})$ vanishes.  We thus have a
$1$-coboundary $\ifam{\eta_a}_{a\in A}\in\cim{1}(\sU,\ssections[r]{\man{E}})$
such that
\begin{equation*}
\eta_b|\nbhd{U}_a\cap\nbhd{U}_b-\eta_a|\nbhd{U}_a\cap\nbhd{U}_b=\xi_{ab},
\qquad a,b\in A.
\end{equation*}
Let $\sigma_a\in\sections[r][\nbhd{U}_a]{\man{B}}$ be given by
$\sigma_a=\sigma_{a0}+\eta_a$ and note that
\begin{equation*}
\sigma_a|\nbhd{U}_a\cap\nbhd{U}_b=(\xi_{a0}+\eta_a)|\nbhd{U}_a\cap\nbhd{U}_b=
(\xi_{b0}+\eta_b)|\nbhd{U}_a\cap\nbhd{U}_b=\sigma_b|\nbhd{U}_a\cap\nbhd{U}_b.
\end{equation*}
Since the sheaf of sections of $\man{B}$ is a sheaf, there exists
$\sigma\in\sections[r]{\man{B}}$ such that $\sigma|\nbhd{U}_a=\sigma_a$\@,
$a\in A$\@.  Thus $\sigma$ is the section we are after.
\end{proof}
\end{proposition}

We have the following two corollaries that will be useful.
\begin{corollary}\label{cor:AsimE}
Let\/ $r\in\{\infty,\omega,\hol\}$\@, and let\/ $\map{\pi}{\man{E}}{\man{M}}$
be a\/ $\C^r$-vector bundle and\/ $\map{\beta}{\man{B}}{\man{M}}$ be a\/
$\C^r$-affine bundle modelled on\/ $\man{E}$\@.  Assume that\/ $\man{M}$ is
Stein when\/ $r=\hol$\@.  Then there exists a\/ $\C^r$-affine bundle
isomorphism\/ $\map{\Psi}{\man{B}}{\man{E}}$ over\/ $\id_{\man{M}}$\@.
\begin{proof}
We first note that, if there exists a $\C^r$-section of
$\map{\beta}{\man{B}}{\man{M}}$\@, then the lemma holds.  Indeed, suppose
that we have a $\C^r$-section $\map{\sigma}{\man{M}}{\man{B}}$\@.  Then one
readily verifies that the mapping
\begin{equation*}
\mapdef{\what{\Psi}}{\man{E}}{\man{B}}{e_x}{\alpha(\sigma(x),e_x)}
\end{equation*}
is a $\C^r$-affine bundle isomorphism, the isomorphism on fibres being that
where $\sigma(x)$ serves as the ``origin'' for the affine space
$\man{B}_x$\@.  We take $\Psi=\what{\Psi}^{-1}$\@.
\end{proof}
\end{corollary}

\begin{corollary}\label{cor:section-exist}
Let\/ $r\in\{\infty,\omega,\hol\}$\@, and let\/ $\map{\pi}{\man{E}}{\man{M}}$
be a\/ $\C^r$-vector bundle and\/ $\map{\beta}{\man{B}}{\man{M}}$ be a\/
$\C^r$-affine bundle modelled on\/ $\man{E}$\@.  Assume that\/ $\man{M}$ is
additionally Stein when\/ $r=\hol$\@.  If\/ $\Xi\in\jet{m}{\man{B}}$\@, then
there exists\/ $\sigma\in\sections[r]{\man{B}}$ satisfying\/
$j_m\sigma(\beta(b))=\Xi$\@.
\begin{proof}
A generalisation of this is proved for vector bundles as Sublemma~2 in the
proof of Theorem~4.5 in~\cite{ADL:20a}\@.  For affine bundles, the result is
then a consequence of Corollary~\ref{cor:AsimE}\@.
\end{proof}
\end{corollary}

\subsection{Jet bundles}

We will consider various sorts of jet bundles in this paper.  We refer
to~\cite{DJS:89} and~\cite[\S{}12]{IK/PWM/JS:93} as useful references.  Here
we shall mainly introduce the notation we use.

Throughout this section, we let $r\in\{\infty,\omega,\hol\}$ and let
$\field=\real$ if $r\in\{\infty,\omega\}$ and let $\field=\complex$ if
$r=\hol$\@.

\subsubsection{Jets of sections of a fibred manifold}

Let $r\in\{\infty,\omega,\hol\}$\@.  We consider a $\C^r$-fibred manifold
$\map{\rho}{\man{X}}{\man{M}}$ (\ie~a surjective submersion) and $\C^r$-local
sections of this manifold.  For $p\in\man{X}$ and $m\in\integernn$\@, we
denote by $\jet[p]{m}{\man{X}}$ the set of $m$-jets of local sections that
take the value $p$ at $x=\rho(p)$\@.  Thus an element of
$\jet[p]{m}{\man{X}}$ is an equivalence class of local sections taking the
value $p$ at $x$ and whose first $m$-derivatives (say, in a fibred chart)
agree at $x$\@.  For a $\C^r$-local section $\sigma$ defined in a
neighbourhood of $x\in\man{M}$\@, we denote by
$j_m\sigma(x)\in\jet[\sigma(x)]{m}{\man{X}}$ the $m$-jet of $\sigma$\@.  We
denote by $\jet{m}{\man{X}}=\disjointunion_{p\in\man{X}}\jet[p]{m}{\man{M}}$
the bundle of $m$-jets of local sections.  Note that
$\jet{0}{\man{X}}\simeq\man{X}$\@.  For $m,l\in\integernn$ with $m\ge l$\@,
we denote by $\map{\rho^m_l}{\jet{m}{\man{X}}}{\jet{l}{\man{X}}}$ the
projection, and we abbreviate
$\map{\rho_m\eqdef\rho\scirc\rho^m_0}{\jet{m}{\man{X}}}{\man{M}}$\@.  This
defines $\jet{m}{\man{X}}$ as a fibred manifold over $\man{M}$\@.

Let us denote by $\vb{\man{X}}=\ker(\tf{\rho})\subset\tb{\man{X}}$ the
vertical bundle of the projection, and let
$\nu=\tbproj{\man{X}}|\vb{\man{X}}$\@.  We claim that
$\map{\rho^m_{m-1}}{\jet{m}{\man{X}}}{\jet{m-1}{\man{X}}}$ is an affine
bundle modelled on
$\rho_{m-1}^*\Symalg*[m]{\ctb{\man{M}}}\otimes(\rho^{m-1}_0)^*\vb{\man{X}}$\@.
The affine structure is defined as follows.  Let $p\in\man{X}$\@, let
$x=\rho(p)$\@, let $\alpha^1,\dots,\alpha^m\in\ctb[x]{\man{M}}$\@, and let
$v\in\vb[p]{\man{X}}$\@.  Let $f^1,\dots,f^m$ be $\C^r$-functions defined
near $x$\@, vanishing at $x$\@, and satisfying $\d{f^j}(x)=\alpha^j$\@,
$j\in\{1,\dots,m\}$\@.  Let $(y,t)\mapsto\omega(y,t)\in\man{X}$ be a mapping
defined near $(x,0)$ and satisfying $\omega(x,t)=p$ for all $t$\@,
$\rho\scirc\omega(y,t)=y$ for all $(y,t)$\@, and
$\tderivatzero{}{t}\omega(x,t)=v$\@.\footnote{One can think of $\omega$ as
being a family of sections that map $x$ to $p$\@.}  Then we define the affine
action by
\begin{equation*}
j_k\sigma(x)+\alpha^1\odot\dots\odot\alpha^m\otimes v=j_k(\sigma_\omega)(x),
\end{equation*}
where $\sigma_\omega$ is the section $\sigma_\omega(x)=\omega(x,f^1(x)\cdots
f^m(x))$\@.  Thus we have the diagram
\begin{equation*}
\xymatrix{{0}\ar[r]&{\rho_{m-1}^*\Symalg*[m]{\ctb{\man{M}}}\otimes
(\rho^{m-1}_0)^*\vb{\man{X}}}\ar[r]&{\jet{m}{\man{X}}}\ar[r]&
{\jet{m-1}{\man{X}}}\ar[r]&{0}}
\end{equation*}
Also note that $\map{\rho^m_0}{\jet{m}{\man{X}}}{\man{X}}$ is an affine
bundle modelled on $\map{\nu_m}{\jet{m}{\vb{\man{X}}}}{\man{X}}$\@.

\subsubsection{Jets of mappings of manifolds}

Let $r\in\{\infty,\omega,\hol\}$ and let $\field\in\{\real,\complex\}$\@, as
appropriate.  We next work with $\C^r$-manifolds $\man{M}$ and $\man{N}$\@,
and jets of $\C^r$-mappings from $\man{M}$ to $\man{N}$\@.  For
$(x,y)\in\man{M}\times\man{N}$\@, we denote by
$\jet[(x,y)]{m}{(\man{M};\man{N})}$ the set of $m$-jets of mappings that map
$x$ to $y$\@.  Thus these are equivalence classes of mappings from $\man{M}$
to $\man{N}$ that map $x$ to $y$ and whose first $m$-derivatives agree (in a
chart, for instance).  For $\Phi\in\mappings[r]{\man{M}}{\man{N}}$\@,
$j_m\Phi(x)\in\jet[(x,\Phi(x))]{m}{(\man{M};\man{N})}$ denotes the $m$-jet of
$\Phi$ at $x$\@.  We denote
\begin{equation*}
\jet{m}{(\man{M};\man{N})}=\bigdisjointunion_{(x,y)\in\man{M}\times\man{N}}
\jet[(x,y)]{m}{(\man{M};\man{N})},
\end{equation*}
which is the bundle of $m$-jets of mappings from $\man{M}$ to $\man{N}$\@.
Note that $\jet{0}{(\man{M};\man{N})}\simeq\man{M}\times\man{N}$\@.  For
$m,l\in\integernn$ with $m\ge l$\@, we
denote
$\map{\rho^m_l}{\jet{m}{(\man{M};\man{N})}}{\jet{l}{(\man{M};\man{N})}}$ as
the natural projection.  We abbreviate
$\map{\rho_m\eqdef\pr_1\scirc\rho^m_0}{\jet{m}{(\man{M};\man{N})}}
{\man{M}}$\@.

In the special case where $\man{N}=\field$\@,~\ie~when we are dealing with
jets of functions, we denote by
$\jetalg[x]{m}{\man{M}}=\jet[(x,0)]{m}{(\man{M};\field)}$ the $m$-jets of
functions that have the value $0$ at $x$\@.  We note that this is an
$\field$-algebra by
\begin{gather*}
j_mf(x)+j_mg(x)=j_m(f+g)(x),\quad a(j_mf(x))=j_m(af)(x),\\
(j_mf(x))\cdot(j_mg(x))=j_m(fg)(x).
\end{gather*}
We can then think of $\jet[(x,y)]{m}{(\man{M};\man{N})}$ as the set of
$\field$-algebra homomorphisms from $\jetalg[y]{m}{\man{N}}$ to
$\jetalg[x]{m}{\man{M}}$  by
\begin{equation*}
j_m\Phi(x)(j_mg(y))=j_m(\Phi^*g)(x),\qquad
j_k\Phi(x)\in\jet[(x,y)]{m}{(\man{M};\man{N})},\
j_kg(y)\in\jetalg[y]{m}{\man{N}}.
\end{equation*}
One can verify that this correspondence is a bijection, and it allows us to
think of jets of mappings in a concrete algebraic context.

To understand some of the structure of jet bundles of mappings, it is
convenient to treat such jets as a special case of the jets of the fibred
manifold $\map{\rho}{\man{M}\times\man{N}}{\man{M}}$\@, with $\rho=\pr_1$\@.
In this case,
$\sections[r]{\man{M}\times\man{N}}\simeq\mappings[r]{\man{M}}{\man{N}}$ by
the observation that $\sigma(x)=(x,\Phi_\sigma(x))$ for
$\sigma\in\sections[r]{\man{M}\times\man{N}}$ and for the associated mapping
$\Phi_\sigma\in\mappings[r]{\man{M}}{\man{N}}$\@.  In like manner, we
identify $\jet{m}{(\man{M}\times\man{N})}$ and $\jet{m}{(\man{M};\man{N})}$
by
\begin{equation*}
j_m\sigma(x)=(x,j_m\Phi_\sigma(x)).
\end{equation*}
We use $\jet{m}{(\man{M};\man{N})}$ to denote the space of $m$-jets in this
setting, although it is sometimes convenient to think of this as
$\jet{m}{(\man{M}\times\man{N})}$\@, and we shall work with both ways of
thinking things, depending on which is most convenient.  We shall use
$\jet{0}{(\man{M};\man{N})}$\@, noting that this is simply identified with
$\man{M}\times\man{N}$\@.  Note that
$\vb{(\man{M}\times\man{N})}\simeq0\oplus\tb{\man{N}}$\@.  An important
specialisation that happens in this case of the jet bundle of a fibred
manifold comes about because the fibred manifold is trivial (and not just
merely trivialisable).  Because of this
$\map{\rho^m_0}{\jet{m}{(\man{M};\man{N})}}{\jet{0}{(\man{M}\times\man{N})}}$
is a \emph{vector bundle}\@, not just an affine bundle.  Indeed, it is
isomorphic to the vector bundle
$\rho_m^*\jetalg{m}{\man{M}}\otimes(\pr_2\scirc\rho^m_0)^*\tb{\man{N}}$\@.
Note, however, that while $\jet{m}{(\man{M};\man{N})}$ is a vector bundle, it
is not the jet bundle of sections of a vector bundle.  Now, for
$(x,y)\in\jet{0}{(\man{M};\man{N})}$ and $m\in\integerp$\@, we define
$\map{\epsilon_m}{\Symalg*[m]{\ctb[x]{\man{M}}}\otimes\tb[y]{\man{N}}}
{\jet[(x,y)]{m}{(\man{M};\man{N})}}$ by
\begin{equation*}
\epsilon_m(\d{f^1}(x)\odot\dots\odot\d{f^m}(x)\otimes Y(y))(j_kg(y))=
j_m(f^1\cdots f^m(\lieder{Y}{g}(y)))(x)
\end{equation*}
where $f^1,\dots,f^m$ are $\C^r$-functions on $\man{M}$ defined near $x$ and
which vanish at $x$\@, and where $Y$ is a $\C^r$-vector field on $\man{N}$
defined near $y$\@.  Note that, in writing this formula, we are defining an
$m$-jet of mappings as an algebra homomorphism from $\jetalg[y]{m}{\man{N}}$
to $\jetalg[x]{m}{\man{M}}$\@.  This then gives rise to the short exact
sequence
\begin{equation}\label{eq:JMN-sequence}
\xymatrix{{0}\ar[r]&{\rho_{m-1}^*\Symalg*[m]{\ctb[x]{\man{M}}}\otimes
(\pr_2\scirc\rho^{m-1}_0)^*\tb{\man{N}}}\ar[r]^(0.7){\epsilon_m}&
{\jet{m}{(\man{M};\man{N})}}\ar[r]^{\rho^m_{m-1}}&
{\jet{m-1}{(\man{M};\man{N})}}\ar[r]&{0}}
\end{equation}

\subsubsection{Jets of sections of a vector bundle}

Let $r\in\{\infty,\omega,\hol\}$ and let $\map{\pi}{\man{E}}{\man{M}}$ be a
$\C^r$-vector bundle.  For $x\in\man{M}$ and $m\in\integernn$\@,
$\jet[x]{m}{\man{E}}$ denotes the set of $m$-jets of sections of $\man{E}$ at
$x$\@.  For a $\C^r$-section $\xi$ defined in some neighbourhood of $x$\@,
$j_k\xi(x)\in\jet[x]{m}{\man{E}}$ denotes the $m$-jet of $\xi$\@.  We denote
by $\jet{m}{\man{E}}=\disjointunion_{x\in\man{M}}\jet[x]{m}{\man{E}}$ the
bundle of $m$-jets.  For $m,l\in\integernn$ with $m\ge l$\@, we denote by
$\map{\pi^m_l}{\jet{m}{\man{E}}}{\jet{l}{\man{E}}}$ the projection.  Note
that $\jet{0}{\man{E}}\simeq\man{E}$\@.  We abbreviate
$\map{\pi_m\eqdef\pi\scirc\pi^m_0}{\jet{m}{\man{E}}}{\man{M}}$\@, and note
that $\jet{m}{\man{E}}$ has the structure of a vector bundle over
$\man{M}$\@, with addition and scalar multiplication defined by
\begin{equation*}
j_m\xi(x)+j_m\eta(x)=j_m(\xi+\eta)(x),\quad a(j_m\xi(x))=j_m(a\xi)(x)
\end{equation*}
for sections $\xi$ and $\eta$ and for $a\in\field$\@.

For $x\in\man{M}$ and $m\in\integerp$\@, define
$\map{\epsilon_m}{\Symalg*[m]{\ctb[x]{\man{M}}}\otimes\man{E}_x}
{\jet[x]{m}{\man{E}}}$ by
\begin{equation}\label{eq:epsilonmJmE}
\epsilon_m(\d{f^1}(x)\odot\dots\odot\d{f^m}(x)\otimes\xi(x))=
j_m(f^1\cdots f^m\xi),
\end{equation}
where $f^1,\dots,f^m$ are locally defined $\C^r$-functions around $x$ that
vanish at $x$\@.  One can easily show that we then have the following short
exact sequence:
\begin{equation}\label{eq:JE-sequence}
\xymatrix{{0}\ar[r]&{\pi_{m-1}^*\Symalg*[m]{\ctb{\man{M}}}\otimes
\pi_{m-1}^*\man{E}}\ar[r]^(0.75){\epsilon_m}&
{\jet{m}{\man{E}}}\ar[r]^(0.45){\pi^m_{m-1}}&{\jet{m-1}{\man{E}}}\ar[r]&{0}}
\end{equation}

\subsubsection{Jet of sections of an affine bundle}

We also will talk about jet bundles of affine bundles.  Thus we let
$r\in\{\infty,\omega,\hol\}$ and let $\map{\pi}{\man{E}}{\man{M}}$ be a
$\C^r$-vector bundle with $\map{\beta}{\man{B}}{\man{M}}$ a $\C^r$-affine
bundle modelled on $\man{E}$\@.  The affine bundle, being a fibred manifold,
is entitled to its jet bundles.  Then we have vector bundles
$\map{\pi_m}{\jet{m}{\man{E}}}{\man{M}}$ and affine bundles
$\map{\beta_m}{\jet{m}{\man{B}}}{\man{M}}$\@, $m\in\integerp$\@.  We let
$\vb{\man{B}}=\ker(\tf{\beta})\subset\tb{\man{B}}$ be the vertical bundle,
which is a vector bundle $\map{\nu}{\vb{\man{B}}}{\man{B}}$ over $\man{B}$\@.
Note that $\vb{\man{B}}\simeq\beta^*\man{E}$\@.  We note that
$\map{\beta^m_{m-1}}{\jet{m}{\man{B}}}{\jet{m-1}{\man{B}}}$ is an affine
bundle modelled according to the short exact sequence
\begin{equation}\label{eq:JB-sequence}
\xymatrix{{0}\ar[r]&
{\beta_{m-1}^*\Symalg*[m]{\ctb{\man{M}}}\otimes\beta_{m-1}^*\man{E}}
\ar[r]^(0.75){\epsilon_m}&{\jet{m}{\man{B}}}\ar[r]^{\beta^m_{m-1}}&
{\jet{m-1}{\man{B}}}\ar[r]&{0}}
\end{equation}
We claim that $\jet{m}{\man{B}}$ is an affine bundle modelled on
$\jet{m}{\man{E}}$\@.  Indeed, if
$\map{\alpha}{\man{E}\times_{\man{M}}\man{B}}{\man{B}}$ is the affine bundle
structure for $\man{B}$\@, then the affine bundle structure for
$\jet{m}{\man{B}}$ is
\begin{equation}\label{eq:jmalpha}
\mapdef{j_m\alpha}{\jet{m}{\man{E}}\times_{\man{M}}\jet{m}{\man{B}}}
{\jet{m}{\man{B}}}{(j_m\xi(x),j_m\sigma(x))}{j_m(\alpha(\xi,\sigma))(x).}
\end{equation}
Given a section $\sigma\in\sections[r]{\man{B}}$\@, we have the $\C^r$-affine
bundle isomorphism
\begin{equation*}
\mapdef{\iota_\sigma}{\man{B}}{\man{E}}{b}{b-\sigma(b)}
\end{equation*}
from Corollary~\ref{cor:AsimE}\@.  In terms of sections, this gives an isomorphism
\begin{equation*}
\mapdef{\hat{\iota}_\sigma}{\sections[r]{\man{B}}}
{\sections[r]{\man{E}}}{\gamma}{\gamma-\sigma}
\end{equation*}
of $\field$-affine spaces.  This then induces an isomorphism
\begin{equation*}
\mapdef{j_m\iota_\sigma}{\jet{m}{\man{B}}}{\jet{m}{\man{E}}}
{j_m\gamma(x)}{j_m\gamma(x)-j_m\sigma(x)}
\end{equation*}
of affine bundles.

\subsection{Functions on vector and affine bundles}\label{subsec:vbabfunctions}

As we have indicated, a key ingredient in our extensions of the well-known
forms of Gelfand duality in differential geometry is the characterisation of
appropriate classes of functions.  In this section we consider linear and
affine functions.

Let $r\in\{\infty,\omega,\hol\}$ and let $\field\in\{\real,\complex\}$ as
appropriate.  Let $\map{\pi}{\man{E}}{\man{M}}$ be a vector bundle of class
$\C^r$ and let $\map{\beta}{\man{B}}{\man{M}}$ be a $\C^r$-affine bundle
modelled on $\man{E}$\@.  Note that
$\map{\beta^*}{\func[r]{\man{M}}}{\func[r]{\man{B}}}$ is an homomorphism of
$\field$-algebras.  This then gives $\func[r]{\man{B}}$ the structure of a
$\func[r]{\man{M}}$-module with multiplication $f\cdot F=(\beta^*f)F$ for
$f\in\func[r]{\man{M}}$ and $F\in\func[r]{\man{B}}$\@.  Now is a good time to
mention that we will be a little bit sloppy and write either of
\begin{equation*}
f\cdot F,\quad fF,\quad\beta^*fF
\end{equation*}
for the same thing, whichever seems to best illustrate what we are doing at
the moment.  Note that we have a short exact sequence of
$\func[r]{\man{M}}$-modules
\begin{equation}\label{eq:funcMfuncE}
\xymatrix{{0}\ar[r]&{\func[r]{\man{M}}}\ar[r]^{\beta^*}&
{\func[y]{\man{B}}}\ar[r]&{\func[r]{\man{B}}/\beta^*\func[r]{\man{M}}}
\ar[r]&{0}}
\end{equation}
Let us now introduce two particular classes of functions.
\begin{definition}\label{def:fibre-linaff}
Let $r\in\{\infty,\omega,\hol\}$\@, let $\map{\pi}{\man{E}}{\man{M}}$ be a
vector bundle of class $\C^r$\@, and let $\map{\beta}{\man{B}}{\man{M}}$ be a
$\C^r$-affine bundle modelled on $\man{E}$\@.
\begin{compactenum}[(i)]
\item A function $F\in\func[r]{\man{E}}$ is \defn{fibre-linear} if, for each
$x\in\man{M}$\@, $F|\man{E}_x$ is a linear function.
\item A function $F\in\func[r]{\man{B}}$ is \defn{fibre-affine} if, for each
$x\in\man{M}$\@, $F|\man{B}_x$ is an affine function.
\end{compactenum}
We denote by $\linfunc[r]{\man{E}}$ (\resp~$\afffunc[r]{\man{B}}$) the set of
$\C^r$-fibre-linear functions on $\man{E}$ (\resp~$\C^r$-fibre-affine
functions on $\man{B}$).\oprocend
\end{definition}

Let us give some elementary properties of the sets of fibre-linear and
fibre-affine functions.
\begin{lemma}\label{lem:linafffuncs}
Let\/ $r\in\{\infty,\omega,\hol\}$\@, let\/ $\map{\pi}{\man{E}}{\man{M}}$ be
a vector bundle of class\/ $\C^r$\@, and let\/
$\map{\beta}{\man{B}}{\man{M}}$ be an affine bundle of class\/ $\C^r$
modelled on\/ $\man{E}$\@.  Then the following statements hold:
\begin{compactenum}[(i)]
\item \label{pl:fibre-affine1} $\linfunc[r]{\man{E}}$ and\/
$\afffunc[r]{\man{B}}$ are submodules of the\/ $\func[r]{\man{M}}$-modules\/
$\func[r]{\man{E}}$ and\/ $\func[r]{\man{B}}$\@, respectively;

\item \label{pl:fibre-affine2} for\/ $F\in\linfunc[r]{\man{E}}$\@, there
exists\/ $\lambda_F\in\sections[r]{\dual{\man{E}}}$ such that
\begin{equation*}
F(e)=\natpair{\lambda_F\scirc\pi(e)}{e},\qquad e\in\man{E},
\end{equation*}
and, moreover, the map\/ $F\mapsto\lambda_F$ is an isomorphism of\/
$\func[r]{\man{M}}$-modules;

\item \label{pl:fibre-affine3} for\/ $F\in\afffunc[r]{\man{B}}$\@, there
exists\/ $\alpha_F\in\sections[r]{\affdual{\man{B}}}$ such that
\begin{equation*}
F(b)=\natpair{\alpha_F\scirc\beta(b)}{b},\qquad b\in\man{B},
\end{equation*}
and, moreover, the map\/ $F\mapsto\alpha_F$ is an isomorphism of\/
$\func[r]{\man{M}}$-modules;

\item \label{pl:fibre-affine4} the short exact sequence~\eqref{eq:funcMfuncE}
induces a short exact sequence
\begin{equation}\label{eq:funcMafffuncE}
\xymatrix{{0}\ar[r]&{\func[r]{\man{M}}}\ar[r]^{\beta^*}&
{\afffunc[r]{\man{B}}}\ar[r]&{\linfunc[r]{\man{E}}}\ar[r]&{0}}
\end{equation}
of\/ $\func[r]{\man{M}}$-modules;

\item \label{pl:fibre-affine5} there is a splitting\/
$\tau\in\Hom_{\func[r]{\man{M}}}(\linfunc[r]{\man{E}};\afffunc[r]{\man{B}})$
of the preceding short exact sequence.
\end{compactenum}
\begin{proof}
\eqref{pl:fibre-affine1} Let $F\in\afffunc[r]{\man{B}}$ and
$f\in\func[r]{\man{M}}$\@.  Then
\begin{equation*}
f\cdot F(e)=(f\scirc\beta(e))F(e),
\end{equation*}
and so $f\cdot F$ is fibre-affine since a scalar multiple of an affine
function is an affine function.  Also, since the pointwise sum of affine
functions is an affine function, we conclude that $\afffunc[r]{\man{B}}$ is
indeed a submodule of $\func[r]{\man{E}}$\@.  Of course, the same sort of
reasoning applies to fibre-linear functions.

\eqref{pl:fibre-affine2} This merely follows by definition of the dual bundle
$\dual{\man{E}}$\@.

\eqref{pl:fibre-affine3} This merely follows by definition of the affine dual
bundle $\affdual{\man{B}}$\@.

\eqref{pl:fibre-affine4} First note that
$\beta^*(\func[r]{\man{M}})\subset\afffunc[r]{\man{B}}$\@.  Indeed, elements
of $\beta^*(\func[r]{\man{M}})$ are constant on fibres of $\man{E}$\@.  Thus
they are affine with zero linear part.  Now the assertion holds since an
element of $\afffunc[r]{\man{B}}/\beta^*\func[r]{\man{M}}$ consists of
fibre-affine functions that differ by a function that is constant on fibres.
Affine functions differing by a constant have the same linear part, and so we
conclude that elements of $\afffunc[r]{\man{B}}/\beta^*\func[r]{\man{M}}$ are
naturally identified with functions that are linear on fibres.  That is,
\begin{equation*}
\afffunc[r]{\man{B}}/\beta^*\func[r]{\man{M}}\simeq\linfunc[r]{\man{E}},
\end{equation*}
as claimed.

\eqref{pl:fibre-affine5} We note that, by Corollary~\ref{cor:AsimE}\@, the bundles
$\man{B}$ and $\man{E}$ are isomorphic as affine bundles.  Thus, if we can
prove this part of the lemma for the affine bundle $\man{E}$\@, it will
follow for the affine bundle $\man{B}$\@.  The result is clear for the affine
bundle $\man{E}$\@, however, since, if $F\in\afffunc[r]{\man{E}}$\@, then
$F(e)=\natpair{\lambda\scirc\pi(e)}{e}+f\scirc\pi(e)$ for some
$\lambda\in\sections[r]{\dual{\man{E}}}$ and $f\in\func[r]{\man{M}}$\@.  Thus
the splitting is obtained by either injecting $\linfunc[r]{\man{E}}$ into
$\afffunc[r]{\man{E}}$ or projecting from $\afffunc[r]{\man{E}}$ to
$\func[r]{\man{M}}$\@.
\end{proof}
\end{lemma}

The lemma ensures that
\begin{equation}\label{eq:afffunc-decomp}
\afffunc[r]{\man{B}}\simeq
\sections[r]{\dual{\man{E}}}\oplus\func[r]{\man{M}},
\end{equation}
the direct sum being of $\func[r]{\man{M}}$-modules, although this
decomposition is not canonical, except in the case that $\man{B}$ is a vector
bundle.  Moreover, an isomorphism~\eqref{eq:afffunc-decomp} is determined by a
choice of section of $\map{\beta}{\man{B}}{\man{M}}$\@.

We close this section by considering functions induced on vector and jet
bundles.
\begin{definition}\label{def:func1formlift}
Let $r\in\{\infty,\omega,\hol\}$\@, as required.  Let
$\map{\pi}{\man{E}}{\man{M}}$ be a $\C^r$-vector bundle and let\/
$\map{\beta}{\man{B}}{\man{M}}$ be a\/ $\C^r$-affine bundle modelled on\/ $\man{E}$\@.
\begin{compactenum}[(i)]
\item \label{pl:lambdaC1} For
$\alpha\in\sections[r]{\affdual{\man{B}}}$\@, the \defn{vertical
evaluation} of $\alpha$ is $\ve{\alpha}\in\afffunc[r]{\man{B}}$ defined
by $\ve{\alpha}(b_x)=\natpair{\alpha(x)}{b_x}$\@.
\item \label{pl:lambdaC2} For $\lambda\in\sections[r]{\dual{\man{E}}}$\@,
the \defn{vertical evaluation} of $\lambda$ is
$\ve{\lambda}\in\linfunc[r]{\man{E}}$ defined by
$\ve{\lambda}(e_x)=\natpair{\lambda(x)}{e_x}$\@.
\item For $f\in\func[r]{\man{M}}$\@, the \defn{horizontal lift} of $f$ is the
function $\hl{f}\in\func[r]{\man{E}}$ defined by $\hl{f}=\pi^*f$\@.\oprocend
\end{compactenum}
\end{definition}

\subsection{Differential operators}

In this section we first work with a general notion of a differential
operator, one that applies to mappings between manifolds.  We then specialise
to affine and linear differential operators.

\subsubsection{Differential operators on fibred
manifolds}\label{subsubsec:diffops}

To motivate our constructions with differential operators, it is convenient
to work in the more-or-less standard setting of fibred manifolds.  In this
setting, we have the following definition.
\begin{definition}
Let $r\in\{\infty,\omega,\hol\}$\@, and let $\map{\rho}{\man{X}}{\man{M}}$
and $\map{\theta}{\man{Y}}{\man{M}}$ be $\C^r$-fibred manifolds.  Let
$\map{\beta}{\man{B}}{\man{M}}$ be a $\C^r$-affine bundle modelled on a
$\C^r$-vector bundle $\map{\pi}{\man{E}}{\man{M}}$\@.
\begin{compactenum}[(i)]
\item A \defn{$\C^r$-differential operator of order $m$} is a $\C^r$-morphism
of fibred manifolds $\map{P}{\jet{m}{\man{X}}}{\man{Y}}$\@.
\item A $\C^r$-differential operator of order $m$\@,
$\map{P}{\jet{m}{\man{X}}}{\man{B}}$ is \defn{fibre-affine} if, for each
$p\in\man{X}$\@, $P|(\rho_0^m)^{-1}(p)$ is an affine mapping with values in
$\man{B}_{\rho(p)}$\@.
\end{compactenum}
By $\DO[r][m]{\man{X};\man{Y}}$ and $\FADO[r][m]{\man{X};\man{B}}$ the spaces
of $\C^r$-differential operators and fibre-affine differential operators of
order $m$\@, respectively.\oprocend
\end{definition}

In terms of usual notions of differential operators, suppose that
$P\in\DO[r][m]{\man{X};\man{Y}}$ and that $\sigma\in\sections[r]{\man{X}}$\@.
Then we define $\what{P}(\sigma)\in\sections[r]{\man{Y}}$ by asking that
$\what{P}(\sigma)(x)=P(j_m\sigma(x))$\@.

We shall have a particular interest in the case of differential operators
with values in the trivial line bundle $\field_{\man{M}}$\@, which we regard
as an affine bundle when we wish to think of fibre-affine differential
operators.  In this case, we abbreviate
\begin{equation*}
\DO[r][m]{\man{X}}=\DO[r][m]{\man{X};\field_{\man{M}}},\quad
\FADO[r][m]{\man{X}}=\FADO[r][m]{\man{X};\field_{\man{M}}}.
\end{equation*}
We note that, if $P\in\FADO[r][m]{\man{X}}$\@, then
\begin{equation*}
P\scirc j_m\sigma(x)=(x,P_0\scirc j_m\sigma(x))
\end{equation*}
for a fibre-affine function $P_0$ on
$\map{\rho^m_0}{\jet{m}{\man{X}}}{\man{X}}$\@.  Thus the set of fibre-affine
functions on $\jet{m}{\man{X}}$\@, as a bundle over $\man{X}$\@, satisfies
\begin{equation*}
\FADO[r][m]{\man{X}}\simeq\sections[r]{\affdual{(\jet{m}{\man{X}})}}.
\end{equation*}
We shall consistently use the symbol $\FADO[r][m]{\man{X}}$ to denote the set
of fibre-affine functions on $\map{\rho^m_0}{\jet{m}{\man{X}}}{\man{X}}$\@.
Since $\map{\rho^m_0}{\jet{m}{\man{X}}}{\man{X}}$ is an affine bundle
modelled on $\map{\nu_m}{\jet{m}{\vb{\man{X}}}}{\man{X}}$\@, we thus have the
following diagram
\begin{equation*}
\xymatrix{{0}\ar[r]&{\func[r]{\man{X}}}\ar[r]&
{\FADO[r][m]{\man{X}}}\ar@{-->}[r]&
{\sections[r]{\dual{(\jet{m}{\vb{\man{X}}})}}}\ar[r]&{0}}
\end{equation*}
where the dashed arrow indicates the projection onto the linear part.

Let us indicate the setting in which we shall use the preceding development.
We shall work with the setting where $\man{X}=\man{M}\times\man{N}$ and where
$\rho=\pr_1$\@.  In this case, we denote by
$\FADO[r][m]{\man{M}\times\man{N}}$ the $\C^r$-fibre-affine differential
operators of order $m$\@.  In this case, since $\jet{m}{(\man{M};\man{N})}$
is a vector bundle over $\jet{0}{(\man{M};\man{N})}$\@, we can also consider
the special case of fibre-affine differential operators that are indeed
fibre-linear.  In such a case, we denote by
$\FLDO[r][m]{\man{M}\times\man{N}}$ the set of sections of the dual bundle of
$\jet{m}{(\man{M};\man{N})}$ over $\jet{0}{(\man{M};\man{N})}$\@.  We thus
have the exact sequence
\begin{equation*}
\xymatrix{{0}\ar[r]&{\func[r]{\man{M}\times\man{N}}}\ar[r]&
{\FADO[r][m]{\man{M}\times\man{N}}}\ar[r]&
{\FLDO[r][m]{\man{M}\times\man{N}}}\ar[r]&{0}}
\end{equation*}

\subsubsection{Affine and linear differential operators}\label{subsubsec:aff-diffops}

The subject of linear differential operators in vector bundles is
classical~\cite[\S10.1]{LIN:96}\@.  We shall need to generalise this to
affine bundles, and we carry out the fairly straightforward generalisation
here.

We begin with the definitions.
\begin{definition}
Let $r\in\{\infty,\omega,\hol\}$\@, and let $\map{\pi}{\man{E}}{\man{M}}$ be
a $\C^r$-vector bundle with $\map{\beta}{\man{B}}{\man{M}}$ an affine bundle
of class $\C^r$ modelled on $\man{E}$\@.  Also let
$\map{\gamma}{\man{A}}{\man{M}}$ be a $\C^r$-affine bundle modelled on the
$\C^r$-vector bundle $\map{\theta}{\man{F}}{\man{M}}$\@.  Let
$m\in\integernn$\@.
\begin{compactenum}[(i)]
\item A \defn{$\C^r$-linear differential operator of order $m$} from
$\man{E}$ to $\man{F}$ is a $\C^r$-vector bundle mapping
$\map{P}{\jet{m}{\man{E}}}{\man{F}}$\@.
\item A \defn{$\C^r$-affine differential operator of order $m$} from
$\man{B}$ to $\man{A}$ is a $\C^r$-affine bundle mapping
$\map{P}{\jet{m}{\man{B}}}{\man{A}}$\@.
\end{compactenum}
By $\LDO[r][m]{\man{B};\man{A}}$ and $\ADO[r][m]{\man{B};\man{A}}$ we
denote the spaces of $\C^r$-linear differential operators and $\C^r$-affine
differential operators of order $m$\@, respectively.\oprocend
\end{definition}

Let us illustrate how the sort of differential operators we define are
operators that differentiate, in the usual sense.  Let
$P\in\ADO[r][m]{\man{B};\man{A}}$ and let $\sigma\in\sections[r]{\man{B}}$\@.
We define $\what{P}(\sigma)\in\sections[r]{\man{A}}$ by asking that
$\what{P}(\sigma)(x)=P(j_m\sigma(x))$\@.  Similar characterisations are
possible for linear differential operators, of course.  We shall have a
particular interest in the case when $\man{A}$ is the trivial line bundle
$\field_{\man{M}}$\@, which we think of as an affine bundle or a vector
bundle as we need.  We shall abbreviate
\begin{equation*}
\ADO[r][m]{\man{B}}=\ADO[r][m]{\man{B};\field_{\man{M}}},\quad
\LDO[r][m]{\man{E}}=\LDO[r][m]{\man{E};\field_{\man{M}}}.
\end{equation*}
Note that, if $P\in\ADO[r][m]{\man{B}}$\@, then
\begin{equation*}
P\scirc j_m\sigma(x)=(x,P_0\scirc j_m\sigma(x))
\end{equation*}
for a fibre-affine function $P_0$ on
$\map{\beta_m}{\jet{m}{\man{B}}}{\man{M}}$\@.  Thus the set of fibre-affine
functions on $\jet{m}{\man{B}}$ satisfies
\begin{equation*}
\ADO[r][m]{\man{B}}\simeq\sections[r]{\affdual{(\jet{m}{\man{B}})}}.
\end{equation*}
We shall adhere to the convention of denoting these fibre-affine functions by
$\ADO[r][m]{\man{B}}$\@, so fixing one of the three possible pieces of
notation.  In a similar fashion, $\LDO[r][m]{\man{E}}$ is to be thought of as
the set of linear functions on $\jet{m}{\man{E}}$\@, or equivalently as the
set of sections of the dual bundle $\dual{\jet{m}{\man{E}}}$\@.  Note that
the above discussions, and Lemma~\ref{lem:linafffuncs}\@, give the following
exact short sequence:
\begin{equation*}
\xymatrix{{0}\ar[r]&{\func[r]{\man{M}}}\ar[r]&
{\ADO[r][m]{\man{B}}}\ar[r]&{\LDO[r][m]{\man{E}}}\ar[r]&{0}}
\end{equation*}

\section{Gelfand duality for manifolds}\label{sec:Membedding}

In this section we overview and give a unified development of the topological
embedding of a $\C^r$-manifold into the weak dual of the space of
$\C^r$-functions for $r\in\{\infty,\omega,\hol\}$\@.  We also consider the
functorial aspects of this embedding,~\ie~how it behaves relative to
morphisms.

\subsection{Algebras and ideals of functions}

Let $r\in\{\infty,\omega,\hol\}$ and let $\field=\complex$ in the case
$r=\hol$ and let $\field=\real$ in the cases $r\in\{\infty,\omega\}$\@.  Let
$\man{M}$ be a manifold of class $\C^r$\@.  Then $\func[r]{\man{M}}$ has
the structure of an $\field$-algebra, where the algebra structure is that
inherited pointwise from the ring structure of $\field$\@:
\begin{equation*}
(f+g)(x)=f(x)+g(x);\quad(af)(x)=a(f(x));\quad(fg)(x)=f(x)g(x),
\end{equation*}
for $f,g\in\func[r]{\man{M}}$ and $a\in\field$\@.

Motivated by the above, let us make some general comments about
$\field$-algebras.  Thus we let $\alg{A}$ be a commutative $\field$-algebra
with unit $1_{\alg{A}}$\@.  Let $\alg{I}\subset\alg{A}$ be an ideal, thinking
of $\alg{A}$ as a mere ring.  We note that $\alg{A}/\alg{I}$ is an
$\field$-algebra with the algebra operations
\begin{equation*}
(r_1+\alg{I})+(r_2+\alg{I})=(r_1+r_2)+\alg{I},\quad
(r_1+\alg{I})(r_2+\alg{I})=r_1r_2+\alg{I},\quad
a\cdot(r+\alg{I})=a\cdot r+\alg{I}.
\end{equation*}
Note that, given an $\field$-algebra $\alg{A}$\@, we have a canonical
injection $\map{\nu_{\alg{A}}}{\field}{\alg{A}}$ given by
$\nu_{\alg{A}}(a)=a\cdot 1_{\alg{A}}$\@.  One easily verifies that
$\nu_{\alg{A}}$ is an homomorphism of $\field$-algebras.  For
$\field$-algebras we denote by $\AHom_{\field}(\alg{A};\alg{B})$ the set of
$\field$-algebra homomorphisms from the $\field$-algebra $\alg{A}$ to the
$\field$-algebra $\alg{B}$\@.  We denote by $\Aut_{\field}(\alg{A})$ the set
of $\field$-algebra isomorphisms of a $\field$-algebra $\alg{A}$\@.  An
$\field$-algebra homomorphism $\phi\in\AHom_{\field}(\alg{A};\alg{B})$
between $\field$-algebras is \defn{unital} if
$\phi(1_{\alg{A}})=1_{\alg{B}}$\@.  Unital $\field$-valued homomorphisms have
the following useful property.
\begin{lemma}
Let\/ $\field\in\{\real,\complex\}$\@.  For a commutative\/
$\field$-algebra\/ $\alg{A}$ with unit and an ideal\/
$\alg{I}\subset\alg{A}$\@, the following statements are equivalent:
\begin{compactenum}
\item \label{pl:real-ideal1} the map\/
$\map{\nu_{\alg{A}/\alg{I}}}{\field}{\alg{A}/\alg{I}}$ is an isomorphism;
\item \label{pl:real-ideal2} there exists a unital homomorphism\/
$\phi\in\AHom_{\field}(\alg{A};\field)$ of\/ $\field$-algebras for which\/
$\alg{I}=\ker(\phi)$\@.
\end{compactenum}
Moreover, an ideal satisfying the two equivalent conditions is maximal.
\begin{proof}
\eqref{pl:real-ideal1}$\implies$\eqref{pl:real-ideal2} Let us define
$\phi\in\AHom_{\field}(\alg{A};\field)$ by
$\phi(f)=\nu_{\alg{A}/\alg{I}}^{-1}\scirc\pi_{\alg{I}}$\@, where
$\map{\pi_{\alg{I}}}{\alg{A}}{\alg{A}/\alg{I}}$ is the canonical projection.
Note that $\nu_{\alg{A}/\alg{I}}(a)=a\cdot(1_{\alg{A}}+\alg{I})$ and so
$\nu_{\alg{A}/\alg{I}}^{-1}(1_{\alg{A}}+\alg{I})=1$\@.  Thus
\begin{equation*}
\phi(1_{\alg{A}})=\nu_{\alg{A}/\alg{I}}^{-1}\scirc\pi_{\alg{I}}(1_{\alg{A}})=
\nu_{\alg{A}/\alg{I}}^{-1}(1_{\alg{A}}+\alg{I})=1,
\end{equation*}
and so $\phi$ is unital.  Moreover, since $\nu_{\alg{A}/\alg{I}}$ is an
isomorphism,
\begin{equation*}
\ker(\phi)=\ker(\nu_{\alg{A}/\alg{I}}^{-1}\scirc\pi_{\alg{I}})=
\ker(\pi_{\alg{I}})=\alg{I}.
\end{equation*}

\eqref{pl:real-ideal2}$\implies$\eqref{pl:real-ideal1} Suppose that
$a\in\ker(\nu_{\alg{A}/\alg{I}})$\@.  Then
\begin{equation*}
0+\alg{I}=\nu_{\alg{A}/\alg{I}}(a)=a\cdot1_{\alg{A}}+\alg{I},
\end{equation*}
which implies that $a\cdot1_{\alg{A}}\in\alg{I}=\ker(\phi)$\@.  Thus 
\begin{equation*}
0=\phi(a\cdot1_{\alg{A}})=a\phi(1_{\alg{A}})=a.
\end{equation*}
Thus $\nu_{\alg{A}/\alg{I}}$ is injective.  By the first isomorphism
theorem~\cite[Theorem~IV.1.7]{TWH:74} the map,
\begin{equation*}
\alg{A}/\ker(\phi)\ni a+\ker(\phi)\mapsto\phi(a)\in\image(\phi)=\field
\end{equation*}
is an isomorphism.  Thus $\alg{A}/\alg{I}=\alg{A}/\ker(\phi)$ is isomorphic
to the one-dimensional $\field$-algebra $\field$\@.  Moreover,
$\nu_{\alg{A}/\alg{I}}$ is thus an injective mapping into a one-dimensional
$\field$-algebra, and so is an isomorphism.

For the last assertion, let $\alg{J}$ be an ideal of $\alg{A}$ such that
$\alg{I}\subset\alg{J}$\@.  If $\alg{I}\subsetneq\alg{J}$ then we must have
$\alg{A}/\alg{J}\subsetneq\alg{A}/\alg{I}$\@.  Since
$\alg{A}/\alg{I}\simeq\field$\@, this means that $\alg{A}/\alg{J}=\{0\}$\@,
and so $\alg{J}=\alg{A}$\@, giving maximality of~$\alg{I}$\@.
\end{proof}
\end{lemma}

For $x\in\man{M}$\@, we have a unital $\field$-algebra homomorphism
\begin{equation*}
\mapdef{\ev_x}{\func[r]{\man{M}}}{\field}{f}{f(x)}
\end{equation*}
called the \defn{evaluation map}\@.  The evaluation map has useful
topological, as well as algebraic structure.  To describe this, first note
that, as an $\field$-algebra homomorphism, $\ev_x$ is $\field$-linear.  Thus
$\ev_x\in\algdual{\func[r]{\man{M}}}$\@.  We shall also see that it is
continuous, and so is an element of $\topdual{\func[r]{\man{M}}}$\@.  We
shall equip $\topdual{\func[r]{\man{M}}}$ with the weak-$*$ topology, that is
the topology defined by the family of seminorms
\begin{equation*}
\mapdef{p_f}{\topdual{\func[r]{\man{M}}}}{\realnn}
{\alpha}{\snorm{\alpha(f)},}
\end{equation*}
for $f\in\func[r]{\man{M}}$\@.

\subsection{The embedding theorem for manifolds}

Our first main result is now the following.
\begin{theorem}\label{the:Membedding}
Let\/ $r\in\{\infty,\omega,\hol\}$ and let\/ $\field\in\{\real,\complex\}$ as
appropriate.  Let $\man{M}$ be a manifold of class $\C^r$\@, Stein when\/
$r=\hol$\@.  Then the mapping
\begin{equation*}
\mapdef{\ev_{\man{M}}}{\man{M}}{\topdual{\func[r]{\man{M}}}}{x}{\ev_x}
\end{equation*}
is an homeomorphism of\/ $\man{M}$ with the set of unital\/ $\field$-algebra
homomorphisms from\/ $\func[r]{\man{M}}$ to\/ $\field$\@, where the latter
has the topology induced by the weak-$*$ topology.
\begin{proof}
For the first part of the proof, we consider $\ev_{\man{M}}$ as taking values
in $\dual{\func[r]{\man{M}}}$\@,~\ie~taking values in the algebraic dual,
which it obviously does.

First we show that $\ev_{\man{M}}$ is injective.  Suppose that
$x_1,x_2\in\man{M}$ are distinct.  We claim that there exists
$f\in\func[r]{\man{M}}$ such that $f(x_1)=0$ and $f(x_2)=1$\@.  For
$r=\infty$\@, one proves this using bump functions.  In case
$r\in\{\omega,\hol\}$ this is still true, although a partition of unity
argument no longer works.  Instead, we use the fact that the hypotheses of
the theorem ensure that there exists a proper $\C^r$-embedding
$\map{\iota_{\man{M}}}{\man{M}}{\field^N}$ for sufficiently large $N$\@; Let
$\map{\Psi}{\field^N}{\field^N}$ be an affine (and so $\C^r$) isomorphism
such that $\Psi\scirc\iota_{\man{M}}(x_1)=\vect{0}$ and
$\Psi\scirc\iota_{\man{M}}(x_2)=\vect{e}_1$\@, where $\vect{e}_j$\@,
$j\in\{1,\dots,N\}$\@, are the standard basis vectors.  For
$j\in\{1,\dots,N\}$\@, let $\map{\pr_j}{\field^N}{\field}$ be the projection
onto the $j$th component, noting that $\pr_j$ is of class $\C^r$\@.  Then the
function $f=\pr_1\scirc\Psi\scirc\iota_{\man{M}}$ has the desired property
that $f(x_1)=0$ and $f(x_2)=1$\@.  Then we have that
$\ev_{x_1}(f)\not=\ev_{x_2}(f)$\@,~\ie~$\ev_{x_1}\not=\ev_{x_2}$\@.  Thus
$\ev_{\man{M}}$ is injective.

Next we show that $\ev_{\man{M}}$ is surjective.  Let
$\map{\psi}{\func[r]{\man{M}}}{\field}$ be a unital $\field$-algebra
homomorphism.  For $f\in\func[r]{\man{M}}$\@, denote
\begin{equation*}
Z(f)=\setdef{x\in M}{f(x)=0}.
\end{equation*}
We then have the following useful lemma.
\begin{prooflemma}
If\/ $f^1,\dots,f^k\in\ker(\psi)$\@, then\/ $\cap_{j=1}^kZ(f^j)\not=\emptyset$\@.
\begin{subproof}
Suppose that $\cap_{j=1}^kZ(f^j)=\emptyset$\@.  Then, for each
$x\in\man{M}$\@, there exists $j\in\{1,\dots,k\}$ such that $f^j(x)\not=0$\@.
Then, if $[g]_x\in\gfunc[r]{x}{\man{M}}$\@, we can write
$[g]_x=([g]_x[1/f^j]_x)[f^j]_x$\@.  Thus we conclude that
$[f^1]_x,\dots,[f^k]_x$ generate $\gfunc[r]{x}{\man{M}}$ for each
$x\in\man{M}$\@.  Define the surjective sheaf morphism
\begin{equation*}
\mapdef{\Psi}{(\sfunc[r]{\man{M}})^k}{\sfunc[r]{\man{M}}}
{([g^1]_x,\dots,[g^k]_x)}{[g^1]_x[f^1]_x+\dots+[g^k]_x[f^k]_x.}
\end{equation*}
Since the sheaves serving as the domain and codomain of $\Psi$ are coherent
in the real analytic and holomorphic cases (by
\cite[Consequence~A.4.2.1]{HG/RR:84}), we can then use (1)~the vanishing of
sheaf cohomology for sheaves over $\sfunc[\infty]{\man{M}}$
(\cite[Proposition~3.11]{ROW:08}\@, along with
\cite[Examples~3.4(d,e)]{ROW:08} and~\cite[Proposition~3.5]{ROW:08}), in the
case $r=\infty$ or (2)~Cartan's Theorem~B in the cases $r\in\{\omega,\hol\}$
to conclude that the map
\begin{equation*}
\mapdef{\Psi_{\man{M}}}{\func[r]{\man{M}}^k}{\func[r]{\man{M}}}
{(g^1,\dots,g^k)}{g^1f^1+\dots+g^kf^k}
\end{equation*}
is surjective.  In particular, there exists $g^1,\dots,g^k\in\func[r]{\man{M}}$ such that
\begin{equation*}
g^1(x)f^1(x)+\dots+g^k(x)f^k(x)=1,\qquad x\in\man{M}.
\end{equation*}
Since $\ker(\psi)$ is an ideal, we conclude that
$\mathsf{1}_{\man{M}}\in\ker(\psi)$\@, contradicting the assumption that
$\psi$ is unital.
\end{subproof}
\end{prooflemma}

We now can complete the proof of surjectivity of $\ev_{\man{M}}$\@.  Let
$\map{\iota_{\man{M}}}{\man{M}}{\field^N}$ be a proper $\C^r$-embedding.  Let
$\chi^1,\dots,\chi^N\in\func[r]{\man{M}}$ be defined by
\begin{equation*}
\iota_{\man{M}}(x)=(\chi^1(x),\dots,\chi^N(x)).
\end{equation*}
Define $f^j\in\func[r]{\man{M}}$ by
$f^j=\chi^j-\psi(\chi^j)\mathsf{1}_{\man{M}}$\@, $j\in\{1,\dots,N\}$\@.
Clearly $f^j\in\ker(\psi)$\@.  By the lemma,
$\cap_{j=1}^NZ(f^j)\not=\emptyset$\@.  If $x\in Z(f^j)$\@, then the
definition of $f^j$ gives $\chi^j(x)=\psi(\chi^j)$\@.  Thus, if
$x\in\cap_{j=1}^NZ(f^j)$\@, then
\begin{equation*}
\iota_{\man{M}}(x)=(\psi(\chi^1),\dots,\psi(\chi^N)).
\end{equation*}
Since $\iota_{\man{M}}$ is injective, $\cap_{j=1}^NZ(f^j)$ is a singleton,
say $\cap_{j=1}^NZ(f^j)=\{x\}$\@.

Now let $f\in\ker(\psi)$\@.  Then
\begin{equation*}
Z(f)\cap\{x\}=Z(f)\cap\left(\bigcap_{j=1}^NZ(f^j)\right).
\end{equation*}
By the lemma, $Z(f)\cap\{x\}\not=\emptyset$\@, and so we must have
$x\in Z(f)$\@.  As this argument is valid for every $f\in\ker(\psi)$\@, we
conclude that, if $f\in\ker(\psi)$\@, then $f(x)=0$\@.  In other words,
$\ker(\psi)\subset\ker(\ev_x)$\@.  Since $\ker(\psi)$ and $\ker(\ev_x)$ are
both maximal ideals, $\ker(\psi)=\ker(\ev_x)$\@.  Let $f\in\func[r]{\man{M}}$
and define $g=f-f(x)\mathsf{1}_{\man{M}}$\@.  Then $\ev_x(g)=0$ and so
$g\in\ker(\ev_x)=\ker(\psi)$\@.  Thus
\begin{equation*}
0=\psi(g)=\psi(f)-f(x)\quad\implies\quad\psi(f)=f(x),
\end{equation*}
\ie~$\psi=\ev_x$\@.

Note that this establishes a bijection between $\man{M}$ and the unital
$\field$-algebra homomorphisms from $\func[r]{\man{M}}$ to $\field$\@.  It
remains to prove the topological assertions of the theorem.

To this end, let us first show that $\ev_{\man{M}}$ is well-defined, in that
$\ev_x$ is continuous for each $x\in\man{M}$\@.  Note that the
$\C^r$-topology is finer than the $\C^0$-topology, so it suffices to show
that $\ev_x$ is continuous in the $\C^0$-topology on $\func[r]{\man{M}}$\@.
The $\C^0$-topology (by definition) is the locally convex topology with the
seminorms
\begin{equation*}
p^0_{\nbhd{K}}(f)=\sup\setdef{\snorm{f(x)}}{x\in\nbhd{K}},\qquad
\nbhd{K}\subset\man{M}\ \textrm{compact}.
\end{equation*}
Let $\nbhd{K}\subset\man{M}$ be compact such that $x\in\nbhd{K}$\@.  Then
\begin{equation*}
\snorm{\ev_x(f)}=\snorm{f(x)}\le p_{\nbhd{K}}^0(f)
\end{equation*}
for $f\in\func[r]{\man{M}}$\@, which gives continuity of $\ev_x$ in the
$\C^0$-topology.

Now we prove that $\ev_{\man{M}}$ is continuous.  Let $x_0\in\man{M}$ and let
$\nbhd{O}\subset\topdual{\func[r]{\man{M}}}$ be a neighbourhood of
$\ev_{x_0}$\@.  Let $f^1,\dots,f^k\in\func[r]{\man{M}}$ and
$r_1,\dots,r_k\in\realp$ be such that
\begin{equation*}
\bigcap_{j=1}^k\setdef{\alpha\in\topdual{\func[r]{\man{M}}}}
{\snorm{\alpha(f^j)-f^j(x_0)}<r_j}\subset\nbhd{O}.
\end{equation*}
Let $\nbhd{U}$ be a neighbourhood of $x_0$ such that
$\snorm{f^j(x)-f^j(x_0)}<r_j$\@, $j\in\{1,\dots,k\}$\@, $x\in\nbhd{U}$\@.
Then, if $x\in\nbhd{U}$\@, $\ev_{\man{M}}(x)\in\nbhd{O}$ and this gives
continuity of $\ev_{\man{M}}$\@.

Finally, we show that $\ev_{\man{M}}$ is an homeomorphism onto its image.  As
above, we let $\map{\iota_{\man{M}}}{\man{M}}{\field^N}$ be a proper
$\C^r$-embedding of $\man{M}$ in $N$-dimensional Euclidean space and denote
by $\chi^1,\dots,\chi^N$ the coordinate functions restricted to $\man{M}$\@,
which are of class $\C^r$\@.  Let $x_0\in\man{M}$ and let $\nbhd{U}$ be a
neighbourhood of $x_0$ in the standard topology of $\man{M}$\@.  Let
$r\in\realp$ be sufficiently small that
\begin{equation*}
\asetdef{x\in\man{M}}{\sum_{j=1}^N\snorm{\chi^j(x)-\chi^j(x_0)}<r}
\subset\nbhd{U}.
\end{equation*}
Note that
\begin{align*}
\iota_{\man{M}}(\man{M})\cap&
\left(\bigcap_{j=1}^N\asetdef{\alpha\in\topdual{\func[r]{\man{M}}}}
{\snorm{\alpha(\chi^j)-\chi^j(x_0)}<r/N}\right)\\
=&\;\asetdef{\iota_{\man{M}}(x)\in\iota_{\man{M}}(\man{M})}
{\snorm{\chi^j(x)-\chi^j(x_0)}<r/N,\ j\in\{1,\dots,N\}}\\
\subset&\;\asetdef{\iota_{\man{M}}(x)\in\iota_{\man{M}}(\man{M})}
{\sum_{j=1}^N\snorm{\chi^j(x)-\chi^j(x_0)}<r}\subset\iota_{\man{M}}(\nbhd{U}).
\end{align*}
This shows that $\nbhd{U}$ is open in the topology induced by
$\ev_{\man{M}}$\@.  That is to say, the topology on $\man{M}$ induced by
$\ev_{\man{M}}$ is finer than the standard topology, which shows that
$\ev_{\man{M}}$ is open onto its image.
\end{proof}
\end{theorem}

We note that the 1--1 correspondence of $\man{M}$ with the unital
$\field$-algebra homomorphisms does not require any topology.  That is to
say, we do not require that $\ev_{\man{M}}(x)=\ev_x$ be continuous, and the
fact that $\ev_{\man{M}}$ is an homeomorphism onto its image is additional to
the 1--1 correspondence.

We claim that the assignment to an object $\man{M}$ in the category of
$\C^r$-manifolds of the object $\func[r]{\man{M}}$ in the category of
$\field$-algebras is injective.  Indeed, suppose that
$\func[r]{\man{M}_1}=\func[r]{\man{M}_2}$\@, equality being as
$\field$-algebras.  Then certainly these algebras have the same collection of
unital $\field$-algebra homomorphisms.  But then this implies that
$\man{M}_1=\man{M}_2$ by the theorem.

\subsection{Mappings of manifolds and homomorphisms of algebras}

An important facet of Gelfand duality is that it assigns not only algebraic
objects to manifolds, but homomorphisms of these algebraic objects to
mappings of manifolds.

The result we state along these lines is the following.
\begin{theorem}\label{the:Mmorphisms}
Let\/ $r\in\{\infty,\omega,\hol\}$ and let\/ $\field\in\{\real,\complex\}$\@,
as appropriate.  Let\/ $\man{M}$ and\/ $\man{N}$ be\/ $\C^r$-manifolds, Stein
if\/ $r=\hol$\@.  Then the following statements hold:
\begin{compactenum}[(i)]
\item \label{pl:Mmorphisms1} if\/ $\Phi\in\mappings[r]{\man{M}}{\man{N}}$\@,
then\/ $\map{\Phi^*}{\func[r]{\man{N}}}{\func[r]{\man{M}}}$ is a continuous\/
$\field$-algebra homomorphism;
\item \label{pl:Mmorphisms2} the mapping\/ $\Phi\mapsto\Phi^*$ is a bijection
from\/ $\mappings[r]{\man{M}}{\man{N}}$ to\/
$\AHom_{\field}(\func[r]{\man{N}};\func[r]{\man{M}})$\@;
\item \label{pl:Mmorphisms3} if\/ $\Phi\in\Diff[r]{\man{M}}$\@, then\/
$\Phi^*$ is a continuous automorphism of\/ $\func[r]{\man{M}}$\@;
\item \label{pl:Mmorphisms4} the mapping\/ $\Phi\mapsto\Phi^*$ is a bijection
from\/ $\Diff[r]{\man{M}}$ to\/ $\Aut_{\field}(\func[r]{\man{M}})$\@.
\end{compactenum}
\begin{proof}
\eqref{pl:Mmorphisms1} Since $\Phi^*$ is clearly $\field$-linear and since
\begin{equation*}
\Phi^*(fg)(x)=(fg)\scirc\Phi(x)=(f\scirc\Phi(x))(g\scirc\Phi(x))=
\Phi^*f(x)\Phi^*g(x)=(\Phi^*f\Phi^*g)(x),
\end{equation*}
for $f,g\in\func[r]{\man{N}}$ and $x\in\man{M}$\@, we conclude that $\Phi^*$
is an $\field$-algebra homomorphism.  Continuity of $\Phi^*$ is proved
in~\cite{ADL:20a}\@; see the remarks at the end of the background discussion
in Section~\ref{sec:intro}\@.

\eqref{pl:Mmorphisms2} First we show that $\Phi\mapsto\Phi^*$ is an injective
mapping.  Suppose that $\Phi_1^*=\Phi_2^*$\@.  Let
$\map{\iota_{\man{N}}}{\man{N}}{\field^N}$ be a $\C^r$-embedding with
coordinate functions $\chi^1,\dots,\chi^N\in\func[r]{\man{N}}$\@.  Then we
have
\begin{align*}
&\Phi_1^*\chi^j(x)=\Phi_2^*\chi^j(x),\qquad
j\in\{1,\dots,N\},\ x\in\man{M},\\
\implies\enspace&\chi^j\scirc\Phi_1(x)=\chi^j\scirc\Phi_2(x),
\qquad j\in\{1,\dots,N\},\ x\in\man{M},\\
\implies\enspace&\iota_{\man{M}}\scirc\Phi_1(x)=\iota_{\man{M}}\scirc\Phi_2(x),
\qquad x\in\man{M},\\
\implies\enspace&\Phi_1(x)=\Phi_2(x),\qquad x\in\man{M},
\end{align*}
as desired.

To show that $\Phi\mapsto\Phi^*$ is surjective, we shall construct a right
inverse of this mapping.  Thus let
$\gamma\in\AHom_{\field}(\func[r]{\man{N}};\func[r]{\man{M}})$ and define
$\map{\Phi_\gamma}{\man{M}}{\man{N}}$ as follows.  Let $x\in\man{M}$ and
$f,g\in\func[r]{\man{N}}$\@, and note that
\begin{equation*}
\ev_x\scirc\gamma(\mathsf{1}_{\man{N}})=\ev_x(\mathsf{1}_{\man{M}})=1
\end{equation*}
and
\begin{equation*}
\ev_x\scirc\gamma(fg)=\ev_x(\gamma(f)\gamma(g))=\gamma(f)(x)\gamma(g)(x)=
(\ev_x\scirc\gamma(f))(\ev_x\scirc\gamma(g)),
\end{equation*}
from which we conclude that
$\map{\ev_x\scirc\gamma}{\func[r]{\man{N}}}{\field}$ is a unital
$\field$-algebra homomorphism.  Thus, by Theorem~\ref{the:Membedding}\@, there exists
$y_x\in\man{N}$ such that $\ev_{y_x}=\ev_x\scirc\gamma$\@.  We define
$\Phi_\gamma(x)=y_x$\@.

We claim that $\Phi_\gamma\in\mappings[r]{\man{M}}{\man{N}}$\@.  Indeed, let
$f\in\func[r]{\man{N}}$ and note that
\begin{equation*}
\Phi_\gamma^*f(x)=f(y_x)=\ev_{y_x}(f)=\ev_x\scirc\gamma(f)=\gamma(f)(x),
\end{equation*}
\ie~$\Phi^*_\gamma f=\gamma(f)\in\func[r]{\man{M}}$\@.  We claim that this
implies that $\Phi_\gamma$ is of class $\C^r$\@.  Indeed, let $x_0\in\man{M}$
and denote $y_0=\Phi_\gamma(x_0)$\@.  Let $(\nbhd{U},\phi)$ be a coordinate
chart for $\man{M}$ about $x_0$ whose coordinate functions we denote by
$\chi^1,\dots,\chi^n$\@.  Let $(\nbhd{V},\psi)$ be a coordinate chart for
$\man{N}$ about $y_0$ whose coordinate functions $\eta^1,\dots,\eta^k$ are
restrictions of globally defined functions of class $\C^r$\@.  This is
possible by Lemma~\ref{lem:global-coordinates}\@.  The mapping
\begin{equation*}
\mapdef{\vect{\eta}}{\man{N}}{\field^k}{y}{(\eta^1(y),\dots,\eta^k(y))}
\end{equation*}
is a diffeomorphism from a neighbourhood $\nbhd{V}'\subset\nbhd{V}$ of $y_0$
to a neighbourhood $\nbhd{W}$ of $\vect{\eta}(y_0)\in\field^k$\@.  Since
$\Phi_\gamma^*\vect{\eta}$ is continuous by hypothesis, there is a
neighbourhood $\nbhd{U}'$ of $x$ such that
$\Phi_\gamma^*\vect{\eta}(\nbhd{U}')\subset\nbhd{W}$\@.  Thus
$\Phi_\gamma(\nbhd{U}')\subset\nbhd{V}'$\@.  Therefore, we can assume without
loss of generality that $\Phi_\gamma(\nbhd{U})\subset\nbhd{V}$\@.  We denote
\begin{equation*}
\mapdef{\vect{\chi}}{\nbhd{U}}{\field^n}{x}{(\chi^1(x),\dots,\chi^n(x)).}
\end{equation*}
Note that the local representative of $\Phi_\gamma$ in the charts
$(\nbhd{U},\phi)$ and $(\nbhd{V},\psi)$ is
\begin{equation*}
\mapdef{\vect{\Phi}_\gamma}{\phi(\nbhd{U})}{\psi(\nbhd{V})}
{\vect{x}}{\vect{\eta}\scirc\Phi_\gamma\scirc\vect{\chi}^{-1}.}
\end{equation*}
Since $\vect{\eta}\scirc\Phi_\gamma$ is of class $\C^r$ (by hypothesis) and
$\vect{\chi}^{-1}$ is of class $\C^r$\@, the local representative of
$\Phi_\gamma$ is of class $\C^r$\@, and this shows that $\Phi_\gamma$ is of
class $\C^r$\@.

Moreover, the equality $\Phi_\gamma^*=\gamma$ proved above is exactly the
statement that the mapping $\gamma\mapsto\Phi_\gamma$ is a right inverse of
the mapping $\Phi\mapsto\Phi^*$\@, and this completes the proof of this part
of the theorem.

\eqref{pl:Mmorphisms3} This follows from part~\eqref{pl:Mmorphisms1} since the
inverse of $\Phi^*$ is $\Phi_*=(\Phi^{-1})^*$ in the case that $\Phi$ is a
diffeomorphism.

\eqref{pl:Mmorphisms4} This follows from part~\eqref{pl:Mmorphisms2}\@, just as
part~\eqref{pl:Mmorphisms3} follows from~\eqref{pl:Mmorphisms1}\@.
\end{proof}
\end{theorem}

\begin{corollary}
Let\/ $r\in\{\infty,\omega,\hol\}$ and let\/ $\field\in\{\real,\complex\}$\@,
as appropriate.  The category of\/ $\C^r$-manifolds is a full subcategory of
the opposite category of the category of\/ $\field$-algebras via the functor
given by\/ $\man{M}\mapsto\func[r]{\man{M}}$ on objects and by\/
$\Phi\mapsto\Phi^*$ on morphisms.
\end{corollary}

\section{Gelfand duality for vector and affine
bundles}\label{sec:Aembedding}

Of course, Theorem~\ref{the:Membedding} applies just as well to the total space of a
vector or affine bundle, at least when one realises that the total space of a
affine bundle over a Stein manifold is a Stein manifold (which we shall do in
Proposition~\ref{prop:Estein}).  Since we already know that there is a 1--1
correspondence between $\C^r$-manifolds and their $\field$-algebras of
$\C^r$-functions, when specialising to manifolds with additional structure
(such as the structure of a vector or affine bundle), we imagine that we
should restrict ourselves to considerations to subsets of the algebras of
$\C^r$-functions.  Specifically, we work with a particular class of
$\field$-linear mappings $\map{\phi}{\afffunc[r]{\man{B}}}{\field}$\@.  We
require these to be compatible with~(1)~the $\field$-algebra structure on
$\func[r]{\man{M}}$ and~(2)~the affine space structure of the fibres of
$\man{B}$\@.  We must also devise suitable morphisms for this algebraic
structure that capture the essential algebraic features of mappings between
affine bundles.

\subsection{Semialgebras}

Let us first contrive a general setting,~\ie~a category, for the class of
functions in which we are interested.
\begin{definition}
Let $\alg{R}$ be a commutative ring with unit.  An
\defn{$\alg{R}$-semialgebra} is a triple $(\alg{M},\alg{A},\kappa)$ for which
\begin{compactenum}[(i)]
\item $\alg{A}$ is an $\alg{R}$-algebra,
\item $\alg{M}$ is an $\alg{A}$-module, and
\item $\kappa\in\Hom_{\alg{A}}(\alg{A};\alg{M})$\@.
\end{compactenum}
We call $\alg{A}$ the \defn{nonlinear part} of the semialgebra and the
$\alg{A}$-module $\alg{M}/\image(\kappa)$ the \defn{linear part} of the
semialgebra, denoted by $L(\alg{M},\alg{A},\kappa)$\@.\oprocend
\end{definition}

\begin{definition}
A \defn{morphism} of $\alg{R}$-semialgebras $(\alg{M}_1,\alg{A}_1,\kappa_1)$
and $(\alg{M}_2,\alg{A}_2,\kappa_2)$ is a pair $(\phi,\phi_0)$ such that
$\phi\in\Hom_{\alg{R}}(\alg{M}_1;\alg{M}_2)$ and
$\phi_0\in\AHom_{\alg{R}}(\alg{A}_1;\alg{A}_2)$ are such that the diagram of $\alg{R}$-modules
\begin{equation*}
\xymatrix{{\alg{M}_1}\ar[r]^{\phi}&{\alg{M}_2}\\
{\alg{A}_1}\ar[u]^{\kappa_1}\ar[r]_{\phi_0}&{\alg{A}_2}\ar[u]_{\kappa_2}}
\end{equation*}
commutes and such that
\begin{equation}\label{eq:semialg-intertwine}
\phi(a_1x_1)=\phi_0(a_1)\phi(x_1),\qquad a_1\in\alg{A}_1,\ x_1\in\alg{M}_1.
\end{equation}
We denote by
\begin{equation*}
\Hom_{\alg{R}}((\alg{M}_1,\alg{A}_1,\kappa_1);
(\alg{M}_2,\alg{A}_2,\kappa_2))
\end{equation*}
the set of $\alg{R}$-semialgebra morphisms.\oprocend
\end{definition}

We leave to the reader the simple exercise of checking that if
$(\alg{M}_1,\alg{A}_1,\kappa_1)$\@, $(\alg{M}_2,\alg{A}_2,\kappa_2)$\@, and
$(\alg{M}_3,\alg{A}_3,\kappa_3)$ are $\alg{R}$-semialgebras, and if
\begin{equation*}
(\phi,\phi_0)\in\Hom_{\alg{R}}((\alg{M}_1,\alg{A}_1,\kappa_1);
(\alg{M}_2,\alg{A}_2,\kappa_2)),\quad
(\psi,\psi_0)\in\Hom_{\alg{R}}((\alg{M}_2,\alg{A}_2,\kappa_2);
(\alg{M}_3,\alg{A}_3,\kappa_3)),
\end{equation*}
then
\begin{equation*}
(\psi\scirc\phi,\psi_0\scirc\phi_0)\in
\Hom_{\alg{R}}((\alg{M}_1,\alg{A}_1,\kappa_1);
(\alg{M}_3,\alg{A}_3,\kappa_3)).
\end{equation*}
Obviously $(\id_{\alg{M}},\id_{\alg{A}})$ is an $\alg{R}$-semialgebra
morphism from $(\alg{M},\alg{A},\kappa)$ to itself, and it has the usual
attributes of an identity morphism.  In short, we have a category of
$\alg{R}$-semialgebras.  We denote by
$\Aut_{\alg{R}}(\alg{M},\alg{A},\kappa)$ the set of $\alg{R}$-semialgebra
isomorphisms of $(\alg{M},\alg{A},\kappa)$\@.

Morphisms of semialgebras induce morphisms on their linear parts.
\begin{lemma}
Let\/ $\alg{R}$ be a ring, and let\/ $(\alg{M}_1,\alg{A}_1,\kappa_1)$
and\/ $(\alg{M}_2,\alg{A}_2,\kappa_2)$ be\/ $\alg{R}$-semialgebras.  If
\begin{equation*}
(\phi,\phi_0)\in\Hom_{\alg{R}}((\alg{M}_1,\alg{A}_1,\kappa_1);
(\alg{M}_2,\alg{A}_2,\kappa_2)),
\end{equation*}
then there is an induced\/ $\alg{R}$-module homomorphism
\begin{equation*}
L(\phi,\phi_0)\in\Hom_{\alg{R}}(L(\alg{M}_1,\alg{A}_1,\kappa_1);
L(\alg{M}_2,\alg{A}_2,\kappa_2))
\end{equation*}
for which the diagram
\begin{equation*}
\xymatrix{{\alg{A}_1}\ar[r]^{\kappa_1}\ar[d]_{\phi_0}&
{\alg{M}_1}\ar[r]\ar[d]_{\phi}&{L(\alg{M}_1,\alg{A}_1,\kappa_1)}
\ar[d]^{L(\phi,\phi_0)}\\{\alg{A}_2}\ar[r]^{\kappa_2}&{\alg{M}_2}\ar[r]&
{L(\alg{M}_2,\alg{A}_2,\kappa_2)}}
\end{equation*}
of\/ $\alg{R}$-modules commutes and for which
\begin{equation*}
L(\phi,\phi_0)(a_1(x_1+\image(\kappa_1)))=
\phi_0(a_1)L(\phi,\phi_0)(x_1+\image(\kappa_1))
\end{equation*}
for\/ $a_1\in\alg{A}_1$\@,\/ $x_1+\image(\kappa_1)\in L(\alg{M}_1,\alg{A}_1,\kappa_1)$\@.
\begin{proof}
Define
\begin{equation*}
L(\phi,\phi_0)(x_1+\image(\kappa_1))=\phi(x_1)+\image(\kappa_2).
\end{equation*}
Let us show that $L(\phi,\phi_0)$ is well-defined.  This is standard.
Indeed, suppose that $x'_1-x_1=\kappa_1(a_1)$ so that
$x_1+\image(\kappa_1)=x'_1+\image(\kappa_1)$\@.  Then
\begin{align*}
\phi(x'_1)=&\;\phi(x_1+\kappa_1(a_1))=\phi(x_1)+\phi\scirc\kappa_1(a_1)\\
=&\;\phi(x_1)+\kappa_2\scirc\phi_0(a_1),
\end{align*}
showing that $L(\phi,\phi_0)$ is indeed well-defined.  We also have
\begin{align*}
L(\phi,\phi_0)(a_1(x_1+\image(\kappa_1)))=&\;
L(\phi,\phi_0)(a_1x_1+\image(\kappa_1))\\
=&\;\phi(a_1x_1)+\image(\kappa_2)\\
=&\;\phi_0(a_1)\phi(x_1)+\image(\kappa_2)\\
=&\;\phi_0(a_1)(\phi(x_1)+\image(\kappa_2))\\
=&\;\phi_0(a_1)L(\phi,\phi_0)(x_1+\image(\kappa_1)),
\end{align*}
as claimed.
\end{proof}
\end{lemma}

Let us consider the examples of semialgebras.
\begin{examples}\label{eg:semialgebra}
\begin{compactenum}
\item \label{enum:alg->semialg} If $\alg{R}$ is a ring and $\alg{A}$ is an
$\alg{R}$-algebra, we can identify $\alg{A}$ it in a natural way with the
$\alg{R}$-semialgebra $(\alg{A},\alg{A},\id_{\alg{A}})$\@.  The linear part
of such semialgebras is the zero module.  One can see that morphisms of
semialgebras of this form are essentially homomorphisms of algebras since
they arise from diagrams like the following:
\begin{equation*}
\xymatrix{{\alg{A}_1}\ar[r]^{\phi_0}&{\alg{A}_2}\\
{\alg{A}_1}\ar[u]^{\id_{\alg{A}_1}}\ar[r]_{\phi_0}&
{\alg{A}_2}\ar[u]_{\id_{\alg{A}_2}}}
\end{equation*}
The intertwining condition~\eqref{eq:semialg-intertwine} then reads
\begin{equation*}
\phi_0(a_1b_1)=\phi_0(a_1)\phi_0(b_1),\qquad a_1,b_1\in\alg{A}_1,
\end{equation*}
\ie~it expresses that the semialgebra morphisms are algebra morphisms.

\item \label{enum:affbun-semialg} Let $r\in\{\infty,\omega,\hol\}$ and let
$\field\in\{\real,\complex\}$\@, as appropriate. If
$\map{\beta}{\alg{B}}{\man{M}}$ is a $\C^r$-affine bundle modelled on the
$\C^r$-vector bundle $\map{\pi}{\man{E}}{\man{M}}$\@, then
$(\afffunc[r]{\man{B}},\func[r]{\man{M}},\beta^*)$ is an
$\field$-semialgebra.  The nonlinear part of the semialgebra is the algebra
$\func[r]{\man{M}}$ of functions and the linear part is the
$\func[r]{\man{M}}$-module $\afffunc[r]{\man{B}}/\func[r]{\man{M}}\simeq\linfunc[r]{\man{E}}$\@.

If $\map{\beta_1}{\man{B}_1}{\man{M}_1}$ and
$\map{\beta_2}{\man{B}_2}{\man{M}_2}$ are $\C^r$-affine bundles and if
$\map{\Phi}{\man{B}_1}{\man{B}_2}$ is a $\C^r$-affine bundle map over
$\map{\Phi_0}{\man{M}_1}{\man{M}_2}$\@, then we claim that
$(\Phi^*,\Phi_0^*)$ is an $\field$-semialgebra morphism from
$(\afffunc[r]{\man{B}_2},\func[r]{\man{M}_2},\beta_2^*)$ to
$(\afffunc[r]{\man{B}_1},\func[r]{\man{M}_1},\beta_1^*)$\@.  The
commutativity of the diagram
\begin{equation*}
\xymatrix{{\man{B}_1}\ar[r]^{\Phi}\ar[d]_{\beta_1}&{\man{B}_2}\ar[d]^{\beta_2}\\
{\man{M}_1}\ar[r]_{\Phi_0}&{\man{M}_2}}
\end{equation*}
gives the commuting of the diagram
\begin{equation*}
\xymatrix{{\afffunc[r]{\man{B}_2}}\ar[r]^{\Phi^*}&{\afffunc[r]{\man{B}_1}}\\
{\func[r]{\man{M}_1}}\ar[u]^{\beta_1^*}\ar[r]_{\Phi_0^*}&
{\func[r]{\man{M}_2}}\ar[u]_{\beta_2^*}}
\end{equation*}
Let us show that $\Phi^*F$ is fibre-affine if $F\in\afffunc[r]{\man{B}_2}$.
Let $x_1\in\man{M}_1$ and note that
\begin{equation*}
\Phi^*F|\man{B}_{1,x_1}=(F|\man{B}_{2,\Phi_0(x_1)})\scirc\Phi|\man{B}_{1,x_1}.
\end{equation*}
Since $F$ is fibre-affine and $\Phi$ is affine, we see that, indeed,
$\Phi^*F$ is fibre-affine.  Finally, if $f_2\in\func[r]{\man{M}_2}$ and
$F_2\in\afffunc[r]{\man{B}_2}$\@, then
\begin{equation*}
\Phi^*(f_2F_2)(b_1)=(f_2F_2)(\Phi(b_1))=
f_2(\Phi_0(\beta(b_1)))F_2*\Phi(b_1),
\end{equation*}
which gives the intertwining condition~\eqref{eq:semialg-intertwine}\@.\oprocend
\end{compactenum}
\end{examples}

The last example makes it evident what we will be going for here; we will
arrive at the category of affine bundles being a full subcategory of the
category of $\field$-semialgebras.  Analogously to our presentation in
Section~\ref{sec:Membedding}\@, for $b\in\man{B}$ we have the $\field$-linear
mapping
\begin{equation*}
\mapdef{\Ev_b}{\afffunc[r]{\man{B}}}{\field}{F}{F(b).}
\end{equation*}
Clearly $\Ev_b$ is unital,~\ie~$\Ev_b(\mathsf{1}_{\man{B}})=1$\@.  We claim
that
\begin{equation*}
(\Ev_b,\ev_{\beta(b)})\in\Hom_{\field}
((\afffunc[r]{\man{B}},\func[r]{\man{M}},\beta^*);
(\field,\field;\id_{\field})).
\end{equation*}
Indeed, the diagram
\begin{equation*}
\xymatrix{{\afffunc[r]{\man{M}}}\ar[r]^(0.6){\Ev_b}&{\field}\ar@{=}[d]\\
{\func[r]{\man{M}}}\ar[u]^{\beta^*}\ar[r]_(0.6){\ev_{\beta(b)}}&{\field}}
\end{equation*}
clearly commutes.  Also,
\begin{equation*}
\Ev_b(fF)=f(\beta(b))F(b)=\ev_{\beta(b)}(f)\Ev_b(F),\qquad
f\in\func[r]{\man{M}},\ F\in\afffunc[r]{\man{B}},
\end{equation*}
and so the intertwining condition~\eqref{eq:semialg-intertwine} also holds.
We shall often grammatically identify $(\Ev_b,\ev_{\beta(b)})$ with $\Ev_b$
for convenience.  Thus, for example, we may say that $(\Ev_b,\ev_{\beta(b)})$
is an element of $\algdual{\afffunc[r]{\man{B}}}$\@, when we mean to apply
this assertion only to $\Ev_b$\@.

We note that, as a particular facet of its definition, we have
$\Ev_b\in\dual{\afffunc[r]{\man{B}}}$\@,~\ie~$\Ev_b$ is a member of the
algebraic dual.  We shall see, moreover, that it is a member of the
topological dual $\topdual{\afffunc[r]{\man{B}}}$\@.  We equip
$\topdual{\afffunc[r]{\man{B}}}$ with the weak-$*$ topology defined by the
family of seminorms
\begin{equation*}
\mapdef{P_F}{\topdual{\afffunc[r]{\man{B}}}}{\realnn}
{\nu}{\snorm{\nu(F)},}
\end{equation*}
for $F\in\afffunc[r]{\man{B}}$\@.  We similarly use the weak-$*$ topology for
$\topdual{\linfunc[r]{\man{E}}}$\@.

\subsection{The embedding theorem for affine bundles}

The following result now gives us the desired topological embedding for
affine bundles.
\begin{theorem}\label{the:Aembedding}
Let\/ $r\in\{\infty,\omega,\hol\}$ and let\/ $\field\in\{\real,\complex\}$ as
appropriate.  Let\/ $\map{\pi}{\man{E}}{\man{M}}$ be a vector bundle of
class\/ $\C^r$\@, assuming that\/ $\man{M}$ is Stein in the case\/
$r=\hol$\@.  Let\/ $\map{\beta}{\man{B}}{\man{M}}$ be an affine bundle
modelled on\/ $\man{E}$\@.  Then the mapping
\begin{equation*}
\mapdef{\Ev_{\man{B}}}{\man{B}}{\topdual{\afffunc[r]{\man{B}}}}{b}{\Ev_b}
\end{equation*}
is an homeomorphism of\/ $\man{B}$ with the set of unital\/
$\field$-semialgebra morphisms from\/
$(\afffunc[r]{\man{B}},\func[r]{\man{M}},\beta^*)$ to\/
$(\field,\field,\id_{\field})$\@, where the latter has the topology induced
by the weak-$*$ topology.  Moreover,\/ $\Ev_{\man{B}}|\man{B}_x$ is an affine
map for each\/ $x\in\man{M}$\@.
\begin{proof}
Let us show that $\Ev_{\man{B}}$ is injective.  Let $b_1,b_2\in\man{B}$ be
distinct.  As in Corollary~\ref{cor:abembedding}\@, let
$\map{\iota_{\man{B}}}{\man{B}}{\field^N\times\field^N}$ be an injective
$\C^r$-affine bundle mapping over a proper $\C^r$-embedding
$\map{\iota_{\man{M}}}{\man{M}}{\field^N}$\@.

Suppose first that $\beta(b_1)\not=\beta(b_2)$\@.  Then, as in the proof of
Theorem~\ref{the:Membedding}\@, there exists $f\in\func[r]{\man{M}}$ such that
$f\scirc\beta(b_1)=0$ and $f\scirc\beta(b_2)=1$\@.  Then
$\beta^*f\in\afffunc[r]{\man{B}}$ is such that $\beta^*f(b_1)=0$ and
$\beta^*f(b_2)=1$\@.  Thus $\Ev_{\man{E}}(b_1)\not=\Ev_{\man{E}}(b_2)$ in this
case.

Now suppose that $\beta(b_1)=\beta(b_2)$\@.  Let
$\map{\Psi}{\field^N}{\field^N}$ be a linear isomorphism (thus of class
$\C^r$) such that $\Psi\scirc\pr_2\scirc\iota_{\man{B}}(b_1)=\vect{e}_1$ and
$\Psi\scirc\pr_2\scirc\iota_{\man{B}}(b_2)=\vect{e}_2$\@, where
$\vect{e}_1,\dots,\vect{e}_N$ are the standard basis vectors.  Denote
\begin{equation*}
\mapdef{\ol{\Psi}}{\field^N\times\field^N}{\field^N\times\field^N}
{(\vect{x},\vect{v})}{(\vect{x},\Psi(\vect{v})).}
\end{equation*}
Note that $\ol{\Psi}$ is a $\C^r$-affine bundle isomorphism over
$\id_{\field^N}$\@.  For $j\in\{1,\dots,N\}$\@, denote
\begin{equation*}
\mapdef{\Pr_j}{\field^N\times\field^N}{\field}
{((x_1,\dots,x_N),(v_1,\dots,v_N))}{v_j.}
\end{equation*}
Note that $\Pr_j$ is of class $\C^r$\@.  Define $F\in\func[r]{\man{B}}$ by
$F=\Pr_1\scirc\ol{\Psi}\scirc\iota_{\man{B}}$\@.  Note that $F$ is
fibre-affine since $\iota_{\man{B}}$ and $\ol{\Psi}$ are injective affine
bundle mappings and since $\Pr_1$ is fibre-affine.  Moreover, $F(b_1)=1$ and
$F(b_2)=0$\@.  Thus we again conclude that
$\Ev_{\man{E}}(b_1)\not=\Ev_{\man{E}}(b_2)$\@.

Next we prove that $\Ev_{\man{B}}$ is surjective.  Let $(\psi,\psi_0)$ be a
unital $\field$-semialgebra morphism from
$(\afffunc[r]{\man{B}},\func[r]{\man{M}},\beta^*)$ to
$(\field,\field,\id_{\field})$\@.
We then have the following useful lemma.
\begin{prooflemma}
If\/ $F^1,\dots,F^k\in\ker(\psi)$\@, then\/ $\cap_{j=1}^kZ(f^j)\not=\emptyset$\@.
\begin{subproof}
Suppose that $\cap_{j=1}^kZ(F^j)=\emptyset$\@.  Let us write
$F^j=\ve{(\lambda^j)}$ for $\lambda^j\in\sections[r]{\affdual{\man{B}}}$\@.
Then, for each $x\in\man{M}$ and for each $b\in\man{B}_x$\@, there exists
$j\in\{1,\dots,k\}$ such that $F^j(b)\not=0$\@.  Thus
$\lambda^1(x),\dots,\lambda^k(x)$ span $\affdual{\man{B}}_x$ and so
$[\lambda^1]_x,\dots,[\lambda^k]_x$ generate
$\gsections[r]{x}{\affdual{\man{B}}}$\@.  Thus we have a surjective sheaf
morphism
\begin{equation*}
\mapdef{\Psi}{(\sfunc[r]{\man{M}})^k}{\ssections[r]{\affdual{\man{B}}}}
{([g^1]_x,\dots,[g^k]_x)}{[g^1]_x[\lambda^1]_x+\dots+[g^k]_x[\lambda^k]_x.}
\end{equation*}
Since the sheaves serving as the domain and codomain of $\Psi$ are coherent
in the real analytic and holomorphic cases (by
\cite[Consequence~A.4.2.1]{HG/RR:84}), we can then use (1)~the vanishing of
sheaf cohomology for sheaves over $\sfunc[\infty]{\man{M}}$
(\cite[Proposition~3.11]{ROW:08}\@, along with
\cite[Examples~3.4(d,e)]{ROW:08} and~\cite[Proposition~3.5]{ROW:08}), in the
case $r=\infty$ or (2)~Cartan's Theorem~B in the cases $r\in\{\omega,\hol\}$
to conclude that the map
\begin{equation*}
\mapdef{\Psi_{\man{B}}}{\func[r]{\man{M}}^k}{\afffunc[r]{\man{B}}}
{(g^1,\dots,g^k)}{\beta^*g^1F^1+\dots+\beta^*g^kF^k}
\end{equation*}
is surjective.  In particular, there exists $g^1,\dots,g^k\in\func[r]{\man{M}}$ such that
\begin{equation*}
\beta^*g^1F^1+\dots+\beta^*g^kF^k=\mathsf{1}_{\man{B}}=
\beta^*\mathsf{1}_{\man{M}}.
\end{equation*}
Thus
\begin{align*}
\psi(\beta^*g^1F^1+\dots+\beta^*g^kF^k)=&\;
\psi(\beta^*g^1F^1)+\dots+\psi(\beta^*g^kF^k)\\
=&\;\psi_0(g^1)\psi(F^1)+\dots+\psi_0(g^k)\psi(F^k)=0.
\end{align*}
This is in contradiction with $\psi$ being unital.
\end{subproof}
\end{prooflemma}

We now can complete the proof of surjectivity of $\Ev_{\man{B}}$\@.  Let
$\map{\iota_{\man{B}}}{\man{B}}{\field^N\times\field^N}$ be an injective
$\C^r$-affine bundle mapping over a proper $\C^r$-embedding
$\map{\iota_{\man{M}}}{\man{M}}{\field^N}$\@.  Let
\begin{equation*}
((\chi^1,\dots,\chi^N),(\mu^1,\dots,\mu^N))
\end{equation*}
be the coordinate functions for this embedding with
$\chi^1,\dots,\chi^N\in\func[r]{\man{M}}$ and
$\mu^1,\dots,\mu^N\in\afffunc[r]{\man{B}}$\@.  Define
$f^j\in\func[r]{\man{M}}$ by
$f^j=\chi^j-\psi_0(\chi^j)\mathsf{1}_{\man{M}}$\@, $j\in\{1,\dots,N\}$\@.
Define $F^j\in\afffunc[r]{\man{B}}$ by
$F^j=\mu^j-\psi(\mu^j)\mathsf{1}_{\man{B}}$\@.  Clearly $f^j\in\ker(\psi_0)$
and $F^j\in\ker(\psi)$\@.  Note that $\beta^*f^j\in\ker(\psi)$ thanks to the
diagram
\begin{equation}\label{eq:psipsi0}
\xymatrix{{\func[r]{\man{M}}}\ar[r]^{\beta^*}\ar[rd]_{\psi_0}&
{\afffunc[r]{\man{B}}}\ar[d]^{\psi}\\&{\field}}
\end{equation}
By the lemma,
\begin{equation*}
\left(\bigcap_{j=1}^NZ(\beta^*f^j)\right)\cap
\left(\bigcap_{j=1}^NZ(F^j)\right)\not=\emptyset.
\end{equation*}
If $b\in Z(\beta^*f^j)$\@, then the definition of $f^j$ gives
$\chi^j(\beta(b)))=\psi_0(\chi^j)$\@.  If $b\in Z(F^j)$\@, then the
definition of $F^j$ gives $\mu^j(b)=\psi(\mu^j)$\@.  Thus, if
\begin{equation*}
b\in\left(\bigcap_{j=1}^NZ(\beta^*f^j)\right)\cap
\left(\bigcap_{j=1}^NZ(F^j)\right),
\end{equation*}
then
\begin{equation*}
\iota_{\man{B}}(b)=((\psi_0(\chi^1),\dots,\psi_0(\chi^N)),
(\psi(\mu^1),\dots,\psi(\mu^N))).
\end{equation*}
Since $\iota_{\man{B}}$ is injective, there exists $b\in\man{B}$ such that
\begin{equation*}
\left(\bigcap_{j=1}^NZ(\beta^*f^j)\right)\cap
\left(\bigcap_{j=1}^NZ(F^j)\right)=\{b\}.
\end{equation*}

Now let $F\in\ker(\psi)$\@.  Then
\begin{equation*}
Z(F)\cap\{b\}=Z(F)\cap\left(\bigcap_{j=1}^NZ(\beta^*f^j)\right)\cap
\left(\bigcap_{j=1}^NZ(\mu^j)\right).
\end{equation*}
By the lemma, $Z(F)\cap\{b\}\not=\emptyset$\@, and so we must have
$b\in Z(F)$\@.  As this argument is valid for every $F\in\ker(\psi)$\@, we
conclude that, if $F\in\ker(\psi)$\@, then $F(b)=0$\@.  In other words,
$\ker(\psi)\subset\ker(\Ev_b)$\@.  In particular, if
$\beta^*f\in\ker(\psi)$\@, then
\begin{equation*}
0=\psi(\beta^*f)=\psi_0(f).
\end{equation*}
Thus, $\ker(\psi_0)\subset\ker(\ev_{\beta(b)})$\@.  We then argue as in the
proof of Theorem~\ref{the:Membedding} that $\psi_0=\ev_{\beta(b)}$\@.  Let us next
show that $\ker(\Ev_b)\subset\ker(\psi)$\@.  By Corollary~\ref{cor:section-exist}\@,
let $\sigma\in\sections[r]{\man{B}}$ be such that $\sigma(\beta(b))=b$\@.
Then, for $F\in\ker(\Ev_b)$\@,
\begin{equation*}
0=F(b)=F\scirc\sigma(\beta(b)).
\end{equation*}
Thus $F\scirc\sigma\in\ker(\ev_{\beta(b)})=\ker(\psi_0)$\@, and so
$F\scirc\sigma\in\ker(\psi\scirc\beta^*)$ by~\eqref{eq:psipsi0}\@.  Thus
\begin{equation*}
0=\psi\scirc\beta^*(F\scirc\sigma)=\psi(F\scirc\sigma\scirc\beta)=\psi(F),
\end{equation*}
\ie~if $F\in\ker(\Ev_b)$\@, then $F\in\ker(\psi)$\@.

Finally, let $F\in\afffunc[r]{\man{B}}$ and define
$G=F-F(b)\mathsf{1}_{\man{B}}$\@.  Then $\Ev_b(G)=0$ and so
$G\in\ker(\Ev_b)=\ker(\psi)$\@.  Thus
\begin{equation*}
0=\psi(G)=\psi(F)-F(b)\quad\implies\quad\psi(F)=F(b),
\end{equation*}
\ie~$\psi=\Ev_b$\@.

Next we show that the restriction of $\Ev_{\man{B}}$ to fibres is affine.
Let $x\in\man{M}$ and define
\begin{equation*}
\map{A_x}{\man{B}_x}{\dual{\afffunc[r]{\man{B}}}}
\end{equation*}
by
\begin{equation*}
\natpair{A_x(b)}{F}=\natpair{\Ev_{\man{B}}(b)}{F}=F(b).
\end{equation*}
Similarly, define
\begin{equation*}
\map{L_x}{\man{E}_x}{\dual{\linfunc[r]{\man{E}}}}
\end{equation*}
by
\begin{equation*}
\natpair{L_x(e)}{F}=\natpair{\Ev_{\man{E}}(e)}{F}=F(e).
\end{equation*}
To show that $A_x$ is affine, we will show that, for $b\in\man{B}_x$ and
$e\in\man{E}_x$\@, $e\mapsto A_x(b+e)-A_x(b)$ is (1)~independent of $b$ and
(2)~linear in $e$\@.  Thus, let $b\in\man{B}_x$ and let $e\in\man{E}_x$\@,
and let $F\in\afffunc[r]{\man{B}}$\@.  We have
\begin{align*}
\natpair{A_x(b+e)}{F}=&\;\natpair{\Ev_{\man{B}}(b+e)}{F}=F(b+e)\\
=&\;F(b)+L(F)(e)=\natpair{A_x(b)}{F}+\natpair{L_x(e)}{L(F)},
\end{align*}
where $L(F)\in\linfunc[r]{\man{E}}$ is the linear part of $F$\@.  Thus
\begin{equation*}
\natpair{A_x(b+e)}{F}-\natpair{A_x(b)}{F}=\natpair{L_x(e)}{L(F)}.
\end{equation*}
This verifies that $\Ev_{\man{B}}|\man{B}_x$ is affine, as asserted.

The above shows that $\Ev_{\man{B}}$ is a bijection between $\man{B}$ and the
unital $\field$-semialgebra morphisms from
$(\afffunc[r]{\man{B}},\func[r]{\man{M}},\beta^*)$ to
$(\field,\field,\id_{\field})$\@.  It remains to prove the topological
assertions of the theorem.

To show that $\Ev_{\man{E}}$ is well-defined, one must show that
$\Ev_b\in\topdual{\afffunc[r]{\man{B}}}$\@.  Note that the
$\C^r$-topology is finer than the $\C^0$-topology, so it suffices to show
that $\Ev_b$ is continuous in the $\C^0$-topology on
$\afffunc[r]{\man{B}}$\@.  Let $\nbhd{K}\subset\man{B}$ be compact such that
$b\in\nbhd{K}$\@.  Then there exists $C\in\realp$ such that
\begin{equation*}
\snorm{\Ev_b(F)}=\snorm{F(b)}\le p_{\nbhd{K}}^0(F)
\end{equation*}
for $F\in\afffunc[r]{\man{B}}$\@, giving the desired continuity.  (Here
$p_{\nbhd{K}}^0$ is a seminorm for the $\C^0$-topology, as described in the
proof of Theorem~\ref{the:Membedding}\@.)

Next we show that $\Ev_{\man{B}}$ is continuous.  Let $b_0\in\man{B}$ and let
$\nbhd{O}\subset\topdual{\afffunc[r]{\man{B}}}$ be a neighbourhood of
$\Ev_{b_0}$\@.  Let $F^1,\dots,F^k\in\afffunc[r]{\man{B}}$ and
$r_1,\dots,r_k\in\realp$ be such that
\begin{equation*}
\setdef{\alpha\in\topdual{\afffunc[r]{\man{B}}}}
{\snorm{\alpha(F^j)-F^j(b_0)}<r_j,\ j\in\{1,\dots,k\}}\subset\nbhd{O}.
\end{equation*}
Let $\nbhd{V}$ be a neighbourhood of $b_0$ such that
$\snorm{F^j(b)-F^j(b_0)}<r_j$\@, $j\in\{1,\dots,k\}$\@, $b\in\nbhd{V}$\@.
Then, if $b\in\nbhd{V}$\@, $\Ev_{\man{B}}(b)\in\nbhd{O}$ and this gives
continuity of $\Ev_{\man{B}}$\@.

Finally, we show that $\Ev_{\man{B}}$ is an homeomorphism onto its image.  As
above, we let $\map{\iota_{\man{B}}}{\man{B}}{\field^N\times\field^N}$ be an
injective $\C^r$-affine bundle mapping over a proper $\C^r$-embedding
$\map{\iota_{\man{M}}}{\man{M}}{\field^N}$\@.  Define functions
$\chi^j\in\func[r]{\man{M}}$ and $\mu^j\in\func[r]{\man{B}}$\@,
$j\in\{1,\dots,N\}$\@, by the requirement that
\begin{equation*}
\iota_{\man{B}}(b)=((\chi^1\scirc\beta(b),\dots,\chi^N\scirc\beta(b)),
(\mu^1(b),\dots,\mu^N(b))).
\end{equation*}
Let $b_0\in\man{B}$ and let $\nbhd{V}$ be a neighbourhood of $b_0$ in the
standard topology of $\man{B}$\@.  Let $r\in\realp$ be sufficiently small
that
\begin{equation*}
\asetdef{b\in\man{B}}{\sum_{j=1}^N(\snorm{\chi^j\scirc\beta(b)-
\chi^j\scirc\beta(b_0)}+\snorm{\mu^j(b)-\mu^j(b_0)})<r}
\subset\nbhd{V}.
\end{equation*}
Note that
\begin{align*}
\iota_{\man{B}}(\man{B})&
\cap\left(\bigcap_{j=1}^N\asetdef{\alpha\in\topdual{\afffunc[r]{\man{B}}}}
{\snorm{\alpha(\beta^*\chi^j)-\chi^j\scirc\beta(b_0)},
\snorm{\alpha(\mu^j)-\mu^j(b_0)},<r/2N}\right)\\
=&\;\asetdef{\iota_{\man{B}}(b)\in\iota_{\man{B}}(\man{B})}
{\snorm{\chi^j\scirc\beta(b)-\chi^j\scirc\beta(b_0)},
\snorm{\mu^j(b)-\mu^j(b_0)}<r/2N,\ j\in\{1,\dots,N\}}\\
\subset&\;\asetdef{\iota_{\man{B}}(b)\in\iota_{\man{B}}(\man{B})}
{\sum_{j=1}^N(\snorm{\chi^j\scirc\beta(b)-
\chi^j\scirc\beta(b_0)}+\snorm{\mu^j(b)-\mu^j(b_0)})<r}
\subset\iota_{\man{B}}(\nbhd{V}).
\end{align*}
This shows that $\nbhd{V}$ is open in the topology induced by
$\Ev_{\man{B}}$\@.  That is to say, the topology on $\man{B}$ induced by
$\Ev_{\man{B}}$ is finer than the standard topology, which shows that
$\Ev_{\man{B}}$ is open onto its image.
\end{proof}
\end{theorem}

The theorem applies, obviously, to the special case of vector bundles.  Note
that, when working with vector bundles, one \emph{cannot} take the class of
functions to be the fibre-linear functions $\linfunc[r]{\man{E}}$\@.  The
theorem holds in this case with
\begin{equation*}
\map{\Ev_{\man{E}}}{\man{E}}{\topdual{\afffunc[r]{\man{E}}},}
\end{equation*}
\ie~the theorem is ``the same'' for vector bundles as for affine bundles.
The reason is that fibre-linear functions are unable to encode any
information about the base manifold $\man{M}$\@, whereas the affine functions
contain functions on $\man{M}$ as a submodule.  The only distinguishing
feature of vector bundles as compared to affine bundles in this setup is
that, for vector bundles, one has a canonical decomposition
\begin{equation*}
\afffunc[r]{\man{E}}=\func[r]{\man{M}}\oplus\linfunc[r]{\man{E}}.
\end{equation*}

We note that the 1--1 correspondence of $\man{B}$ with the unital,
$\field$-semialgebra morphisms does not require any topology.  That is to
say, we do not require that $\Ev_{\man{B}}(b)=\Ev_b$ be continuous, and the
fact that $\Ev_{\man{B}}$ is an homeomorphism onto its image is additional to
the 1--1 correspondence.

We claim that the assignment to an object $\man{B}$ in the category of
$\C^r$-affine bundles of the object
$(\afffunc[r]{\man{B}},\func[r]{\man{M}},\beta^*)$ in the category of
$\field$-semialgebras is injective.  Indeed, suppose that
\begin{equation*}
(\afffunc[r]{\man{B}_1},\func[r]{\man{M}_1},\beta_1^*)=
(\afffunc[r]{\man{B}_2},\func[r]{\man{M}_2},\beta_2^*),
\end{equation*}
with equality being as $\field$-semialgebras.  Then
$\func[r]{\man{M}_1}=\func[r]{\man{M}_2}$ as $\field$-algebras, and so
$\man{M}_1=\man{M}_2$\@, as we argued following Theorem~\ref{the:Membedding}\@.
Also, the set of unital $\field$-linear mappings of $\afffunc[r]{\man{B}_1}$
must agree with those of $\afffunc[r]{\man{B}_2}$\@, whence
$\man{B}_1=\man{B}_2$ by virtue of the theorem.  Since
$\beta_1^*=\beta_2^*$\@, $\beta_1=\beta_2$\@.  Also since
$\beta_1^*=\beta_2^*$\@, the linear parts of the semialgebras must agree, but
the linear parts are $\linfunc[r]{\man{E}_1}$ and $\linfunc[r]{\man{E}_2}$\@,
where $\man{E}_1$ and $\man{E}_2$ are the model vector bundles.  This,
however, means that $\dual{\man{E}_1}=\dual{\man{E}_2}$\@, whence
$\man{E}_1=\man{E}_2$\@.

\subsection{Mappings of vector and affine bundles and homomorphisms of
algebras of affine functions}

To complete the story of Gelfand duality for affine bundles, we need to
ensure that morphisms in the category of affine bundles and in the category
of $\field$-semialgebras behave as do those for the category of manifolds.
To this end, we have the following analogue of Theorem~\ref{the:Mmorphisms}\@.
\begin{theorem}\label{the:Amorphisms}
Let\/ $r\in\{\infty,\omega,\hol\}$ and let\/ $\field\in\{\real,\complex\}$\@,
as appropriate.  Let\/ $\map{\alpha}{\man{A}}{\man{M}}$ and\/
$\map{\beta}{\man{B}}{\man{N}}$ be\/ $\C^r$-affine bundles.  When\/
$r=\hol$\@, we assume the base manifolds\/ $\man{M}$ and\/ $\man{N}$ are
Stein.  Then the following statements hold:
\begin{compactenum}[(i)]
\item \label{pl:Amorphisms1} if\/
$(\Phi,\Phi_0)\in\abmappings[r]{\man{A}}{\man{B}}$\@, then
\begin{equation*}
(\Phi^*,\Phi^*_0)\in
\Hom_{\field}((\afffunc[r]{\man{A}},\func[r]{\man{M}},\alpha^*);
(\afffunc[r]{\man{B}},\func[r]{\man{N}},\beta^*)),
\end{equation*}
and\/ $\Phi$ and\/ $\Phi_0$ are continuous;
\item \label{pl:Amorphisms2} the mapping\/ $(\Phi,\Phi_0)\mapsto(\Phi^*,\Phi_0^*)$ is a bijection
from\/ $\abmappings[r]{\man{A}}{\man{B}}$ to
\begin{equation*}
\Hom_{\field}((\afffunc[r]{\man{A}},\func[r]{\man{M}},\alpha^*);
(\afffunc[r]{\man{B}},\func[r]{\man{N}},\beta^*));
\end{equation*}
\item \label{pl:Amorphisms3} taking\/ $\man{N}=\man{M}$\@, if\/
$(\Phi,\id_{\man{M}})$ is an isomorphism of the affine bundles\/ $\man{A}$
and\/ $\man{B}$\@, then
\begin{equation*}
(\Phi^*,\id_{\man{M}}^*)\in
\Hom_{\field}((\afffunc[r]{\man{A}},\func[r]{\man{M}},\alpha^*);
(\afffunc[r]{\man{B}},\func[r]{\man{M}},\beta^*)),
\end{equation*}
is an $\field$-semialgebra isomorphism and\/ $\Phi^*$ is continuous;
\item \label{pl:Amorphisms4} the mapping\/
$(\Phi,\id_{\man{M}})\mapsto(\Phi^*,\id_{\man{M}}^*)$ is a bijection from\/
the set of affine bundle isomorphisms of\/ $\man{A}$ to the set\/
$\Aut_{\field}(\afffunc[r]{\man{B}},\func[r]{\man{M}},\beta^*)$\@.
\end{compactenum}
\begin{proof}
\eqref{pl:Amorphisms1} We verified in Example~\enumdblref{eg:semialgebra}{enum:affbun-semialg} that $(\Phi^*,\Phi_0^*)$ is a morphism of $\field$-semialgebras.  Continuity of $\Phi^*$ and $\Phi_0^*$ follows from~\cite[Theorem~9.3]{ADL:20a} in the manner explained in the corresponding part of the proof of Theorem~\ref{the:Mmorphisms}\@.

\eqref{pl:Amorphisms2} First we show that
$(\Phi,\Phi_0)\mapsto(\Phi^*,\Phi_0^*)$ is injective.  Suppose that
$(\Phi_1^*,\Phi_{1,0}^*)=(\Phi_2^*,\Phi_{2,0}^*)$\@.  Let
$\map{\iota_{\man{B}}}{\man{B}}{\field^N\times\field^N}$ be an injective
$\C^r$-affine bundle mapping with coordinate functions
\begin{equation*}
((\chi^1,\dots,\chi^N),(\mu^1,\dots,\mu^N))
\end{equation*}
satisfying $\chi^1,\dots,\chi^N\in\func[r]{\man{N}}$ and
$\mu^1,\dots,\mu^N\in\afffunc[r]{\man{B}}$\@.  As in the proof of
Theorem~\ref{the:Mmorphisms}\@, we have $\Phi_{1,0}=\Phi_{2,0}$\@.  We also have
\begin{align*}
&\Phi_1^*\mu^j(b)=\Phi_2^*\mu^j(b),
\qquad j\in\{1,\dots,N\},\ b\in\man{B},\\
\implies\enspace&\mu^j\scirc\Phi_1(b)=\mu^j\scirc\Phi_2(b),
\qquad j\in\{1,\dots,N\},\ b\in\man{B},\\
\implies\enspace&\iota_{\man{B}}\scirc\Phi_1(b)=\iota_{\man{B}}\scirc\Phi_2(b),
\qquad b\in\man{B},\\
\implies\enspace&\Phi_1(b)=\Phi_2(b),\qquad b\in\man{B},
\end{align*}
which, combined with the fact that $\Phi_{1,0}=\Phi_{2,0}$\@, gives the
desired conclusion.

To show that $(\Phi,\Phi_0)\mapsto(\Phi^*,\Phi_0^*)$ is surjective, we shall
construct a right inverse of this mapping.  Thus let
\begin{equation*}
(\gamma,\gamma_0)\in
\Hom_{\field}((\afffunc[r]{\man{B}},\func[r]{\man{N}},\beta^*);
(\afffunc[r]{\man{A}},\func[r]{\man{M}},\mu^*)),
\end{equation*}
and let $\Phi_{\gamma_0}\in\mappings[r]{\man{M}}{\man{N}}$ be as in the proof
of Theorem~\ref{the:Mmorphisms}\@; thus $\gamma_0=\Phi_{\gamma_0}^*$\@.  Now define
$\map{\Phi_\gamma}{\man{A}}{\man{B}}$ as follows.  Let $a\in\man{A}_x$ and
$F\in\afffunc[r]{\man{B}}$\@, and note that
\begin{equation*}
\Ev_a\scirc\gamma(\mathsf{1}_{\man{B}})=
\ev_x\scirc\gamma_0(\mathsf{1}_{\man{N}})=\ev_x(\mathsf{1}_{\man{M}})=1.
\end{equation*}
Also, for $g\in\func[r]{\man{N}}$\@,
\begin{equation*}
\Ev_a\scirc\gamma(\beta^*g)=\Ev_a\scirc\gamma_0(g)=
\ev_x\scirc\Phi^*_{\gamma_0}(g),
\end{equation*}
giving the commuting of the diagram
\begin{equation*}
\xymatrix{{\func[r]{\man{N}}}\ar[r]^{\beta^*}
\ar[rd]_{\ev_x\scirc\Phi^*_{\gamma_0}}&{\afffunc[r]{\man{B}}}
\ar[d]^{\Ev_a\scirc\gamma}\\&{\field}}
\end{equation*}
For $g\in\func[r]{\man{N}}$ and $G\in\afffunc[r]{\man{B}}$\@, we have
\begin{equation*}
\Ev_a\scirc\gamma(gG)=\Ev_a(\gamma_0(g)\gamma(G))=
(\ev_x\scirc\Phi^*_{\gamma_0}(g))(\Ev_a\scirc\gamma(G)),
\end{equation*}
which is a verification of the intertwining
condition~\eqref{eq:semialg-intertwine}\@.  From all of this, we conclude that
$(\Ev_a\scirc\gamma,\ev_x\scirc\Phi^*_{\gamma_0})$ is a unital
$\field$-semialgebra morphism from
$(\afffunc[r]{\man{B}},\func[r]{\man{N}},\beta^*)$ to
$(\field,\field,\id_{\field})$\@.  Thus, by Theorem~\ref{the:Aembedding}\@, there
exists $b_a\in\man{N}$ such that $\Ev_{b_a}=\Ev_a\scirc\gamma$\@.  We define
$\Phi_\gamma(a)=b_a$\@.

We claim that
$(\Phi_\gamma,\Phi_{\gamma_0})\in\vbmappings[r]{\man{A}}{\man{B}}$\@.  As we
have already seen in the proof of Theorem~\ref{the:Mmorphisms}\@,
$\gamma=\Phi_{\gamma_0}^*$ and
$\Phi_{\gamma_0}\in\mappings[r]{\man{M}}{\man{N}}$\@.  Now let
$F\in\afffunc[r]{\man{B}}$ and note that
\begin{equation*}
\Phi_\gamma^*F(a)=F(b_a)=\Ev_{b_a}(F)=\Ev_a\scirc\gamma(F)=\gamma(F)(a),
\end{equation*}
\ie~$\Phi^*_\gamma F=\gamma(F)\in\afffunc[r]{\man{B}}$\@.  We claim that this
implies that $\Phi_\gamma$ is of class $\C^r$\@.  Indeed, let $a_0\in\man{A}$
and denote $b_0=\Phi_\gamma(a_0)$\@.  Let $(\nbhd{U},\phi)$ be an affine
bundle chart for $\man{B}$ about $\alpha(a_0)$ whose coordinate functions we
denote by
\begin{equation*}
((\chi^1,\dots,\chi^n),(\mu^1,\dots,\mu^k)).
\end{equation*}
Let $(\nbhd{V},\psi)$ be an affine bundle chart for $\man{B}$ about
$\beta(b_0)$ whose coordinate functions
\begin{equation*}
((\eta^1,\dots,\eta^m),(\nu^1,\dots,\nu^l))
\end{equation*}
are restrictions of globally defined functions of class $\C^r$\@
(fibre-affine functions, in the case of $\nu^1,\dots,\nu^l)$\@.  This is
possible by Corollary~\ref{cor:global-abcoordinates}\@.  The mapping
\begin{equation*}
\mapdef{\vect{\eta}\times\vect{\nu}}{\man{B}}
{\field^m\times\field^l}{b}{((\eta^1(\beta(b)),\dots,\eta^m(\beta(b))),
(\nu^1(b),\dots,\nu^k(b)))}
\end{equation*}
is a fibre-affine isomorphism onto its image from an affine bundle coordinate
neighbourhood $\nbhd{V}'\subset\nbhd{V}$ of $b_0$ to a neighbourhood
$\nbhd{W}\times\field^l$ of
$\vect{\eta}(b_0)\times\{\vect{0}\}\in\field^m\times\field^l$\@.  Since
$\Phi_{\gamma_0}^*\vect{\eta}$ and $\Phi_\gamma^*\vect{\nu}$ are continuous
by hypothesis, there is an affine bundle coordinate neighbourhood $\nbhd{U}'$
of $a_0$ such that
\begin{equation*}
\Phi_\gamma(\vect{\eta}\times\vect{\nu})(\nbhd{U}')\subset
\nbhd{W}\times\field^l.
\end{equation*}
Thus $\Phi_\gamma(\nbhd{U}')\subset\nbhd{V}'$\@.  Therefore, we can assume
without loss of generality that $\Phi_\gamma(\nbhd{U})\subset\nbhd{V}$\@.  We
denote
\begin{equation*}
\mapdef{\vect{\chi}\times\vect{\alpha}}{\nbhd{U}}{\field^n\times\field^k}
{x}{((\chi^1(\alpha(a)),\dots,\chi^n(\alpha(a))),
(\mu^1(a),\dots,\mu^k(a))).}
\end{equation*}
Note that the local representative of $\Phi_\gamma$ in the charts
$(\nbhd{U},\phi)$ and $(\nbhd{V},\psi)$ is
\begin{equation*}
\mapdef{\vect{\Phi}_\gamma}{\phi(\nbhd{U})}{\psi(\nbhd{V})}
{\vect{x}}{(\vect{\eta}\times\vect{\nu})\scirc\Phi_\gamma\scirc
(\vect{\chi}\times\vect{\alpha})^{-1}.}
\end{equation*}
Since $(\vect{\eta}\times\vect{\nu})\scirc\Phi_\gamma$ is of class $\C^r$ (by
hypothesis) and $(\vect{\chi}\times\vect{\alpha})^{-1}$ is of class $\C^r$\@,
the local representative of $\Phi_\gamma$ is of class $\C^r$\@, and this
shows that $\Phi_\gamma$ is of class $\C^r$\@.

Moreover, the equality $\Phi_\gamma^*=\gamma$ proved above is exactly the
statement that the mapping $\gamma\mapsto\Phi_\gamma$ is a right inverse of
the mapping $\Phi\mapsto\Phi^*$\@, and this completes the proof of this part
of the theorem.

\eqref{pl:Amorphisms3} This follows from part~\eqref{pl:Amorphisms1} since the
inverse of $\Phi^*$ is $\Phi_*=(\Phi^{-1})^*$ in the case that $\Phi$ is an
affine bundle isomorphism.

\eqref{pl:Amorphisms4} This follows from part~\eqref{pl:Amorphisms2}\@, just as
part~\eqref{pl:Amorphisms3} follows from~\eqref{pl:Amorphisms1}\@.
\end{proof}
\end{theorem}

\begin{corollary}
Let\/ $r\in\{\infty,\omega,\hol\}$ and let\/ $\field\in\{\real,\complex\}$\@,
as appropriate.  The category of\/ $\C^r$-affine bundles is a full
subcategory of the opposite category of the category of\/
$\field$-semialgebras via the functor given by
\begin{equation*}
(\map{\beta}{\man{B}}{\man{M}})\mapsto
(\afffunc[r]{\man{B}},\func[r]{\man{M}},\beta^*)
\end{equation*}
on objects and by\/ $(\Phi,\Phi_0)\mapsto(\Phi^*,\Phi_0^*)$ on morphisms.
\end{corollary}

\section{Gelfand duality for jet bundles}\label{sec:Jembedding}

As an application of the results of Section~\ref{sec:Aembedding}\@, we shall embed
jet bundles into duals of function spaces in a way that respects the
structure of jet bundles.

\subsection{Embedding jet bundles of sections of affine bundles}

In this section we refine the development in the preceding section for affine
bundles to jet bundles of affine bundles.  Thus we suppose that
$r\in\{\infty,\omega,\hol\}$ and that $\map{\beta}{\man{B}}{\man{M}}$ is a
$\C^r$-affine bundle modelled on the $\C^r$-vector bundle
$\map{\pi}{\man{E}}{\man{M}}$\@.  For $m\in\integernn$\@, we have the
$\C^r$-affine bundle $\map{\beta_m}{\jet{m}{\man{B}}}{\man{M}}$\@.  As we saw
in Section~\ref{subsubsec:aff-diffops}\@, the set of fibre-affine functions on this
latter affine bundle is naturally identified with the set
$\ADO[r][m]{\man{B}}$ of $\C^r$-affine differential operators with values in
$\field$\@.  Thus we have the following short exact sequence of
$\func[r]{\man{M}}$-modules:
\begin{equation*}
\xymatrix{{0}\ar[r]&{\func[r]{\man{M}}}\ar[r]^(0.43){\beta_m^*}&
{\ADO[r][m]{\man{B}}}\ar[r]&{\LDO[r][m]{\man{E}}}\ar[r]&{0}}
\end{equation*}
This puts us squarely in the setting of Section~\ref{sec:Aembedding}\@.  That is to
say, we can consider unital $\field$-semialgebra morphisms
\begin{equation*}
(\psi,\psi_0)\in\Hom_{\field}
((\ADO[r][m]{\man{B}},\func[r]{\man{M}},\beta_m^*),
(\field,\field,\id_{\field}))
\end{equation*}
as a subset of the topological dual $\topdual{\ADO[r][m]{\man{B}}}$\@.

Thus we have the following result.
\begin{theorem}\label{the:JkEembedding}
Let\/ $r\in\{\infty,\omega,\hol\}$ and let\/ $\map{\beta}{\man{B}}{\man{M}}$
be a\/ $\C^r$-affine bundle modelled on the\/ $\C^r$-vector bundle\/
$\map{\pi}{\man{E}}{\man{M}}$\@, assuming that $\man{M}$ is Stein in case\/
$r=\hol$\@.  Let\/ $m\in\integernn$\@.  Then the mapping
\begin{equation*}
\mapdef{\Ev^m_{\man{B}}}{\jet{m}{\man{B}}}{\topdual{\ADO[r][m]{\man{B}}}}
{j_m\sigma(x)}{\ev_{j_m\sigma(x)}}
\end{equation*}
is an homeomorphism of\/ $\jet{m}{\man{B}}$ with the set of unital\/
$\field$-semialgebra morphisms from\/
$(\ADO[r][m]{\man{B}},\func[r]{\man{M}},\beta_m^*)$ to\/
$(\field,\field,\id_{\field})$\@, where the latter has the topology induced
by the weak-$*$-topology.  Moreover,\/ $\Ev^m_{\man{B}}|\jet{m}{\man{B}}_x$
is an affine map for each\/ $x\in\man{M}$\@.
\end{theorem}

We note that there are many other affine bundles in this setting arising
from~\eqref{eq:JE-sequence} and~\eqref{eq:JB-sequence}\@.  This will give rise
to corresponding embeddings, and we leave to the reader the chore of
developing the notation required to state the results.  As a hint, we note
that the fibre-affine functions in this setting will be homogeneous
differential operators.

\subsection{Embedding jet bundles of mappings}\label{subsec:JmMNembed}

Next we consider a suitable embedding of the jet bundle of mappings between manifolds.  Here we let $r\in\{\infty,\omega,\hol\}$ and let $\man{M}$ and $\man{N}$ be $\C^r$-manifolds.  Recall that $\jet{0}{(\man{M};\man{N})}=\man{M}\times\man{N}$\@.  For $m\in\integernn$\@, we have the affine bundle $\map{\rho^m_0}{\jet{m}{(\man{M};\man{N})}}{\jet{0}{(\man{M};\man{N})}}$ modelled on the vector bundle $\jet{m}{(\vb{\jet{0}{(\man{M};\man{N})}})}$\@.  As we saw in Section~\ref{subsubsec:diffops}\@, the set of $\C^r$-fibre-affine functions on this affine bundle is identified with the set $\FADO[r][m]{\man{M}\times\man{N}}$ of fibre-affine differential operators of order $m$ with values in $\field$\@.  Therefore, we have the following short exact sequence of $\func[r]{\man{M}}$-modules:
\begin{equation*}
\xymatrix{{0}\ar[r]&{\func[r]{\jet{0}{(\man{M};\man{N})}}} \ar[r]^(0.47){(\rho^m_0)^*}&{\FADO[r][m]{\man{M}\times\man{N}}}\ar[r]& {\sections[r]{\dual{(\jet{m}{\vb{\jet{0}{(\man{M};\man{N})}}})}}}\ar[r]&{0}} \end{equation*}
Again, we are in the setting of Section~\ref{sec:Aembedding}\@, and so we can consider unital $\field$-semialgebra morphisms
\begin{equation*}
(\psi,\psi_0)\in\Hom_{\field} ((\FADO[r][m]{\man{M}\times\man{N}},\func[r]{\jet{0}{(\man{M};\man{N})}}, (\rho^m_0)^*),(\field,\field,\id_{\field}))
\end{equation*}
as a subset of the topological dual $\topdual{\FADO[r][m]{\man{M}\times\man{N}}}$\@.  The embedding result one then has is the following.
\begin{theorem}\label{the:JkMNembedding}
Let\/ $r\in\{\infty,\omega,\hol\}$\@, and let\/ $\man{M}$ and\/ $\man{N}$ be\/ $\C^r$-manifolds, assuming them to be Stein when\/ $r=\hol$\@.  Let\/ $m\in\integernn$\@.  Then the mapping
\begin{equation*}
\mapdef{\Ev^m_{\man{M}\times\man{N}}}{\jet{m}{(\man{M};\man{N})}}
{\topdual{\FADO[r][m]{\man{M}\times\man{N}}}} {j_m\Phi(x)}{\ev_{j_m\Phi(x)}}
\end{equation*}
is an homeomorphism of\/ $\jet{m}{(\man{M};\man{N})}$ with the set of
unital\/ $\field$-semialgebra morphisms from\/
$(\FADO[r][m]{\man{M}\times\man{N}},\func[r]{\man{M}\times\man{N}},
(\rho^m_0)^*)$ to\/ $(\field,\field,\id_{\field})$\@, where the latter has
the topology induced by the weak-$*$ topology.  Moreover,\/
$\Ev^m_{\man{M}\times\man{N}}|(\rho^m_0)^{-1}(x,y)$ is an affine map for
each\/ $x,y\in\man{M}\times\man{N}$\@.
\end{theorem}

We note that there are many other affine bundles in this setting arising
from~\eqref{eq:JMN-sequence}\@.  This will give rise to corresponding
embeddings, and we leave to the reader the pleasure of developing the
notation required to state the results.  As with our hint above for vector
bundles, we comment that the fibre-affine functions arising in this case will
be homogeneous differential operators.

\section{Smooth, real analytic, and holomorphic versions of the
Serre\textendash{}Swan Theorem}

In this section, to wrap up our collection of interconnected results, we
prove a version of the Serre\textendash{}Swan Theorem for vector bundles in
the three regularity categories with which we are working in this paper.

Our proof relies on embedding the total space of a vector bundle in a
suitable Euclidean space.  In the smooth and real analytic cases, this
follows without problem from the embedding theorems in these cases (see, for
example, the proof of Lemma~\ref{lem:global-coordinates} for references).  In the
case of vector bundles over Stein manifolds, that the total space is itself
is Stein is required to use the corresponding embedding theorem.  This is
well-known to be true, and is typically attributed to Serre, and without
reference as near as we can tell.  Related problems are discussed
in~\cite[\S4.21]{FF:11}\@.  In any case, let us prove here the result we
need, since we make substantial use of the corollaries that follow it.
\begin{proposition}\label{prop:Estein}
If\/ $\map{\pi}{\man{E}}{\man{M}}$ is an holomorphic vector bundle over a
Stein manifold, then\/ $\man{E}$ is a Stein manifold.
\begin{proof}
We will show the following three things:
\begin{compactenum}
\item for each $e\in\man{E}$\@, there exists an holomorphic chart for
$\man{E}$ about $e$ whose coordinate functions are globally defined
holomorphic functions;
\item holomorphic functions on $\man{E}$ separate points,~\ie~if
$e_1,e_2\in\man{E}$ are distinct, then there exists
$f\in\func[\hol]{\man{E}}$ such that $f(e_1)\not=f(e_2)$\@;
\item $\man{E}$ is holomorphically convex,~\ie~if $L\subset\man{E}$ is
compact, then the set
\begin{equation*}
\hcohull(L)\eqdef\setdef{e\in\man{E}}{\snorm{f(e)}\le p_L^0(f),\
f\in\func[\hol]{\man{E}}}
\end{equation*}
is compact.
\end{compactenum}
These suffice to show that $\man{E}$ is a Stein manifold by any of the
various definitions.

Let us prove these in order.

Let $e\in\man{E}$\@.  Let $z=\pi(e)$ and let $(\nbhd{U},\phi)$ be an
holomorphic chart for $\man{M}$ about $z$ whose coordinate functions are
globally defined holomorphic functions; this is possible since $\man{M}$ is
Stein.  Let $\ifam{\alpha^1,\dots,\alpha^m}$ be a basis for
$\dual{\man{E}}_z$\@.  By Cartan's Theorem~A, let
$\sigma^1,\dots,\sigma^m\in\sections[\hol]{\dual{\man{E}}}$ be such that
$\sigma^j(z)=\alpha^j$\@, $j\in\{1,\dots,m\}$\@.  Shrink $\nbhd{U}$ so that
$\ifam{\sigma^1(z'),\dots,\sigma^m(z')}$ is a basis for $\dual{\man{E}}_{z'}$
for $z'\in\nbhd{U}$\@.  Then, if $n$ is the dimension of $\man{M}$\@, define
a chart map for $\pi^{-1}(\nbhd{U})$ by
\begin{equation*}
\mapdef{\Phi}{\pi^{-1}(\nbhd{U})}{\complex^n\times\complex^m}
{e_{z'}}{(\phi(z'),(\natpair{\sigma^1(z')}{e_{z'}},\dots,
\natpair{\sigma^m(z')}{e_{z'}})).}
\end{equation*}
This is the required holomorphic chart for $\man{E}$ in a neighbourhood of
$e$ whose coordinate functions are globally defined holomorphic functions.

Let $e_1,e_2\in\man{E}$ be distinct.  If $\pi(e_1)\not=\pi(e_2)$\@, then let
$f\in\func[\hol]{\man{E}}$ be such that
$f\scirc\pi(e_1)\not=f\scirc\pi(e_2)$\@, this being possible since $\man{M}$
is Stein.  Then $\pi^*f(e_1)\not=\pi^*f(e_2)$ and so $\pi^*f$ separates $e_1$
and $e_2$\@.  Now suppose that $\pi(e_1)=\pi(e_2)=z$\@.  Suppose that
$e_1\not=0$\@, without loss of generality.  Let $\alpha\in\dual{\man{E}}_z$
be such that $\alpha(e_1)=1$ and $\alpha(e_2)=0$\@. By Cartan's Theorem~A,
let $\sigma\in\sections[\hol]{\dual{\man{E}}}$ be such that
$\sigma(z)=\alpha$\@.  Define $f\in\func[\hol]{\man{E}}$ by
$f(e)=\natpair{\sigma\scirc\pi(e)}{e}$\@.  Since $f(e_1)=1$ and $f(e_2)=0$\@,
$f$ separates $e_1$ and $e_2$\@.

We shall show that, if $\ifam{e_j}_{j\in\integerp}$ is a sequence in
$\man{E}$ with no accumulation point, there exists $F\in\func[\hol]{\man{E}}$
such that $\limsup_{j\to\infty}\snorm{F(e_j)}=\infty$\@.  First of all, if
the sequence $\ifam{\pi(e_j)}_{j\in\integerp}$ has no accumulation point,
then, since $\man{M}$ is holomorphically convex, there exists
$f\in\func[\hol]{\man{M}}$ such that
$\limsup_{j\to\infty}\snorm{f\scirc\pi(e_j)}=\infty$\@, and since
$\pi^*f\in\func[\hol]{\man{E}}$ this gives the desired conclusion in this
case.  So suppose that $\ifam{\pi(e_j)}_{j\in\integerp}$ has an accumulation
point, and let us pass to a subsequence in order to obtain the assumption
that $\ifam{e_j}_{j\in\integerp}$ has no accumulation point and that that
$\lim_{j\to\infty}\pi(e_j)=x\in\man{M}$\@.  Choose a local trivialisation
$\map{\Phi}{\man{E}|\nbhd{U}}{\nbhd{U}\times\complex^m}$\@, where $\nbhd{U}$
is a neighbourhood of $x$\@.  Let us write
\begin{equation*}
\Phi(e)=(\pi(e),(g^1(e),\dots,g^m(e))),
\end{equation*}
where $g^1,\dots,g^m\in\func[\hol]{\man{E}|\nbhd{U}}$ are linear on
fibres.  Thus, if we define
\begin{equation*}
\alpha^j(e)=(\pi(e),g^j(e)),\qquad j\in\{1,\dots,\},
\end{equation*}
then $\alpha^1,\dots,\alpha^m\in\sections[\hol]{\dual{\man{E}}|\nbhd{U}}$\@.
By Cartan's Theorem~A, there exists
$\sigma_1,\dots,\sigma_k\in\sections[\hol]{\dual{\man{E}}}$ and
$f^1,\dots,f^k\in\func[\hol]{\nbhd{U}}$ such that
\begin{equation*}
\alpha^j=f^1\cdot\sigma_1+\dots+f^k\sigma_k,\qquad j\in\{1,\dots,m\},
\end{equation*}
possibly after shrinking $\nbhd{U}$\@.  By hypothesis, the sequence
$\ifam{\pr_2\scirc\Phi(e_j)}_{j\in\integerp}$ in $\complex^m$ does not have
an accumulation point.  Therefore, we must have that, for some
$a\in\{1,\dots,m\}$\@, the sequence $\ifam{g^a(e_j)}_{j\in\integerp}$ has no
accumulation point.  Since the function $g^a$ is linear on fibres,
$\limsup_{j\to\infty}\snorm{g^a(e_j)}=\infty$\@.  Therefore, for some
$b\in\{1,\dots,k\}$\@, the we must have
$\limsup_{j\to\infty}\snorm{\pr_2\scirc\sigma_a(e_j)}=\infty$\@, furnishing
us with the desired conclusion.
\end{proof}
\end{proposition}

Using the preceding, we can now prove a vector bundle version of the various
embedding theorems.
\begin{proposition}\label{prop:vbembedding}
If\/ $r\in\{\infty,\omega,\hol\}$ and if\/ $\map{\pi}{\man{E}}{\man{M}}$ is
a\/ $\C^r$-vector bundle, with\/ $\man{M}$ Stein when\/ $r=\hol$\@, then
there exist\/ $N\in\integernn$ and an injective\/ $\C^r$-vector bundle
mapping\/ $\map{\iota_{\man{E}}}{\man{E}}{\field^N\times\field^N}$ over a
proper\/ $\C^r$-embedding\/ $\map{\iota_{\man{M}}}{\man{M}}{\field^N}$\@.
\begin{proof}
By~\cite[Lemma~19]{HW:36} in the case of
$r=\infty$\@,~\cite[Theorem~3]{HG:58} in the case of $r=\omega$\@,
and~\cite{RR:54} and Proposition~\ref{prop:Estein} in the holomorphic case, there exists
a proper $\C^r$-embedding $\Psi$ of $\man{E}$ in $\field^N$ for some
$N\in\integerp$\@.  There is then an induced proper $\C^r$-embedding
$\iota_{\man{M}}$ of $\man{M}$ in $\field^N$ by restricting $\Psi$ to the
zero section of $\man{E}$\@.  Let us take the subbundle $\hat{\man{E}}$ of
$\tb{\field^N}|\iota_{\man{M}}(\man{M})$ whose fibre at
$\iota_{\man{M}}(x)\in\iota_{\man{M}}(\man{M})$ is
\begin{equation*}
\hat{\man{E}}_{\iota_{\man{M}}(x)}=\tf[0_x]{\Psi}(\vb[0_x]{\man{E}}).
\end{equation*}
Now recall that $\man{E}\simeq\zeta^*\vb{\man{E}}$\@, where
$\map{\zeta}{\man{M}}{\man{E}}$ is the zero
section~\cite[\cf][page~55]{IK/PWM/JS:93}\@.  Let us abbreviate
$\iota_{\man{E}}=\tf{\Psi}|\zeta^*\vb{\man{E}}$\@.  We then have the
following diagram
\begin{equation*}
\xymatrix{{\man{E}\simeq\zeta^*\vb{\man{E}}}
\ar[d]_{\pi}\ar[r]^{\iota_{\man{E}}}&
{\field^N\times\field^N}\ar[d]^{\pr_2}\\
{\man{M}}\ar[r]_{\iota_{\man{M}}}&{\field^N}}
\end{equation*}
describing an injective mapping $\iota_{\man{E}}$ of $\C^r$-vector bundles
over the proper $\C^r$-embedding $\iota_{\man{M}}$\@, with the image of
$\iota_{\man{E}}$ being $\hat{\man{E}}$\@.  This is the assertion of the
lemma.
\end{proof}
\end{proposition}

Combining this result with Corollary~\ref{cor:AsimE} we have the following.
\begin{corollary}\label{cor:abembedding}
Let\/ $r\in\{\infty,\omega,\hol\}$\@, let\/ $\map{\pi}{\man{E}}{\man{M}}$ be
a\/ $\C^r$-vector bundle, with\/ $\man{M}$ Stein when\/ $r=\hol$\@, and let\/
$\map{\beta}{\man{B}}{\man{M}}$ be a\/ $\C^r$-affine bundle modelled on\/
$\man{E}$\@.  Then there exist\/ $N\in\integernn$ and an injective\/
$\C^r$-affine bundle mapping\/
$\map{\iota_{\man{B}}}{\man{B}}{\field^N\times\field^N}$ over a proper\/
$\C^r$-embedding\/ $\map{\iota_{\man{M}}}{\man{M}}{\field^N}$\@.
\end{corollary}

We also get the following analogue of Lemma~\ref{lem:global-coordinates}\@.
\begin{corollary}\label{cor:global-abcoordinates}
Let\/ $r\in\{\infty,\omega,\hol,\}$ and let\/
$\field\in\{\real,\complex\}$\@, as appropriate.  Let\/
$\map{\beta}{\man{B}}{\man{M}}$ be a\/ $\C^r$-affine bundle modelled on the\/
$\C^r$-vector bundle\/ $\map{\pi}{\man{E}}{\man{M}}$\@, and suppose that\/
$\man{M}$ is Stein when\/ $r=\hol$\@.  Then, for any\/ $x\in\man{M}$\@,
there exist an affine bundle chart\/ $(\nbhd{V},\psi)$ for\/ $\man{B}$ and a
vector bundle chart\/ $(\nbhd{U},\phi)$ for\/ $\man{E}$ whose
coordinate functions
\begin{equation*}
((\chi^1,\dots,\chi^n),(\alpha^1,\dots,\alpha^k))
\end{equation*}
and
\begin{equation*}
((\chi^1,\dots,\chi^n),(\nu^1,\dots,\nu^k))
\end{equation*}
satisfy the following:
\begin{compactenum}[(i)]
\item $\chi^1,\dots,\chi^n$ are restrictions to\/ $\beta(\nbhd{V})$ of
globally defined\/ $\C^r$-functions;
\item $\alpha^1,\dots,\alpha^k$ are restrictions to\/ $\nbhd{V}$ of globally
defined\/ $\C^r$-fibre affine functions;
\item $\nu^1,\dots,\nu^k$ are restrictions to\/ $\nbhd{U}$ of globally
defined\/ $\C^r$-fibre linear functions.
\end{compactenum}
\begin{proof}
The simple idea of the proof of Lemma~\ref{lem:global-coordinates} is easily
adapted to this situation.
\end{proof}
\end{corollary}

We can now state the following variant of the Serre\textendash{}Swan Theorem.
\begin{theorem}\label{the:swan}
Let\/ $r\in\{\infty,\omega,\textup{hol}\}$ and let\/ $\field=\real$ if\/
$r\in\{\infty,\omega\}$ and let\/ $\field=\complex$ if\/ $r=\textup{hol}$\@.
Let\/ $\man{M}$ be a manifold of class\/ $\C^r$\@.  If\/ $r=\hol$ assume
that\/ $\man{M}$ is Stein.  The following statements hold:
\begin{compactenum}[(i)]
\item \label{pl:gen-swan1} if\/ $\map{\pi}{\man{E}}{\man{M}}$ is a vector
bundle of class\/ $\C^r$\@, then $\sections[r]{\man{E}}$ is a finitely
generated projective module over\/ $\func[r]{\man{M}}$\@; that is to say,\/
$\sections[r]{\man{E}}$ is a direct summand of a finitely generated free
module over\/ $\func[r]{\man{M}}$\@;
\item \label{pl:gen-swan2} if\/ $\sM$ is a finitely generated projective
module over\/ $\func[r]{\man{M}}$\@, then\/ $\sM$ is isomorphic to the
module\/ $\sections[r]{\man{E}}$ of\/ $\C^r$-sections of a\/
$\C^r$-generalised subbundle of\/ $\man{E}$\@.
\end{compactenum}
\begin{proof}
\eqref{pl:gen-swan1} By Proposition~\ref{prop:vbembedding}\@, let
$\map{\iota_{\man{E}}}{\man{E}}{\field^N\times^N}$ be an injective
$\C^r$-vector bundle mapping over a proper $\C^r$-embedding
$\map{\iota_{\man{M}}}{\man{M}}{\real^N}$\@.  Thus we have $\man{E}$ as
isomorphic to a subbundle of the trivial bundle
$\field^N_{\man{M}}\eqdef\man{M}\times\field^N$\@.  Let
$\inprod{\cdot}{\cdot}$ be the standard (Hermitian, if $\field=\complex$)
inner product on $\field^N$ which we think of as a vector bundle metric on
$\field^N_{\man{M}}$\@.  Define $\man{G}_x$ to be the orthogonal complement
to $\man{E}_x$\@, noting that $\man{G}$ is then a $\C^r$-subbundle of
$\field^N_{\man{M}}$ and that $\field^N_{\man{M}}=\man{E}\oplus\man{G}$\@.
Let $\map{\pi_1}{\field^N_{\man{M}}}{\man{E}}$ and
$\map{\pi_2}{\field^N_{\man{M}}}{\man{G}}$ be the projections, thought of as
vector bundle morphisms.  Note that $\sections[r]{\field^N_{\man{M}}}$ is
isomorphic, as a $\func[r]{\man{M}}$-module, to $\func[r]{\man{M}}^N$\@.
Moreover, the map from $\sections[r]{\field^N_{\man{M}}}$ to
$\sections[r]{\man{E}}\oplus\sections[r]{\man{G}}$ given by
\begin{equation*}
\vect{\xi}\mapsto(\pi_1\scirc\vect{\xi})\oplus(\pi_2\scirc\vect{\xi})
\end{equation*}
can be directly verified to be an isomorphism of
$\func[r]{\man{M}}$-modules.  In particular, $\sections[r]{\man{E}}$ is a
summand of the free, finitely generated $\func[r]{\man{M}}$-module $\sections[r]{\field^N_{\man{M}}}$\@.

\eqref{pl:gen-swan2} By definition, there exists a module $\sN$ over
$\func[r]{\man{M}}$ such that
\begin{equation*}
\sM\oplus\sN\simeq\underbrace{\func[r]{\man{M}}\oplus\dots\oplus
\func[r]{\man{M}}}_{N\ \textrm{factors}}.
\end{equation*}
The direct sum on the right is naturally isomorphic to the set of
$\C^r$-sections of the trivial vector bundle
$\field^N_{\man{M}}=\man{M}\times\field^N$\@.  Thus we can write
$\sM\oplus\sN=\sections[r]{\field^N_{\man{M}}}$\@.  For $a\in\{1,2\}$\@, let
$\map{\Pi_a}{\sections[r]{\field^N_{\man{M}}}}
{\sections[r]{\field^N_{\man{M}}}}$ be the projection onto the $a$th factor.
As per~\cite[\S6]{EN:67} (essentially), associated with $\Pi_a$ is a vector
bundle map $\map{\pi_a}{\field^N_{\man{M}}}{\field^N_{\man{M}}}$\@.  Since
$\Pi_a\scirc\Pi_a=\Pi_a$ (by virtue of $\Pi_a$ being a projection),
$\pi_a\scirc\pi_a=\pi_a$\@.  To show that $\sM$ is the set of sections of a
vector subbundle of $\field^N_{\man{M}}$ it suffices to show that $\pi_1$ has
constant rank.  One can easily show that $x\mapsto\rank(\pi_{a,x})$ is lower
semicontinuous for $a\in\{1,2\}$\@.  However, since
$\rank(\pi_{1,x})+\rank(\pi_{2,x})=N$ for all $x\in\man{M}$\@, if
$x\mapsto\rank(\pi_{1,x})$ is lower semicontinuous at $x_0$\@, then
$x\mapsto\rank(\pi_{2,x})$ is upper semicontinuous at $x_0$\@.  Thus we
conclude that both of these functions must be continuous at $x_0$\@.  Since
$x\mapsto\rank(\pi_{1,x})$ is integer-valued, it must therefore be constant.
\end{proof}
\end{theorem}

\printbibliography[heading=bibintoc]

\end{document}